\documentclass[onefignum,onetabnum]{siamart250211}

\usepackage{amsmath}
\usepackage{amssymb}
\usepackage{hyperref}
\usepackage[arrow, matrix, curve]{xy}
\usepackage{bm}
\usepackage{yhmath}
\usepackage{color}

\usepackage{hyperref}

\usepackage{amsmath,amssymb,graphicx,textcomp}
\usepackage{mathtools}
\usepackage{cleveref} % clever referencing
\usepackage{verbatim}
\usepackage{placeins}
\usepackage{enumerate}
\usepackage{enumitem}
\usepackage{subcaption}

\DeclareMathOperator{\gammas}{\gamma_{\mathrm{s}}}
\DeclareMathOperator{\gammad}{\gamma_{\mathrm{D}}}

\newcommand{\erfc}{\operatorname{erfc}}
\newcommand{\erf}{\operatorname{erf}}
\newcommand{\Pe}{\operatorname{Pe}}
\newcommand{\nex}{\mathbf{x}}
\newcommand{\ney}{\mathbf{y}}

\newcommand{\nev}{\mathbf{v}}

\newcommand{\neF}{\mathbf{F}}

\newcommand{\neV}{\mathbf{V}}
\newcommand{\neK}{\mathbf{K}}

\newcommand{\dy}{\mathrm{d}\ney}
\newcommand{\dsx}{\mathrm{d}s_\nex}
\newcommand{\dsy}{\mathrm{d}s_\ney}
\newcommand{\jump}[1]{[ \mkern-2.5mu  [ #1 ] \mkern-2.5mu ]}

\newcommand{\average}[1]{\left\{\!\!\left\{ #1\right\} \!\! \right\}}

\newcommand{\Sw}{\hat{\mathrm{S}}_{\omega}}

\newcommand{\Dw}{\hat{\mathrm{D}}_{\omega}}

\newcommand{\Vw}{\hat{\mathrm{V}}_{\omega}}

\newcommand{\Kw}{\hat{\mathrm{K}}_{\omega}}

\newcommand{\W}{\hat{\mathrm{W}}}
\newcommand{\Ilong}{\mathcal{I}^{\text{long}}_{\omega}}
\newcommand{\Ising}{\mathcal{I}^{\text{sing}}_{\omega}}
\newcommand{\Idiff}{\mathcal{I}^{\text{diff}}_{\omega}}

\newcommand{\Gdiff}{G^{\text{diff}}}

\newcommand{\Ilongz}{\mathcal{I}^{\text{long}}_{0}}
\newcommand{\Isingz}{\mathcal{I}^{\text{sing}}_{0}}

\usepackage[normalem]{ulem}
\makeatletter
\newcommand\hide@visible[1]{%
  \bgroup\fboxsep=.3ex\colorbox{Gray}{begin hide}%
  #1\colorbox{Gray}{end hide}\egroup%
}
\newcommand\hide@hidden[1]{%
  \bgroup\fboxsep=.3ex\colorbox{Gray}{hidden text}%
}
\newcommand\hide@invisible[1]{}
\newcommand\makevisible{\let\hide\hide@visible}
\newcommand\makehidden{\let\hide\hide@hidden}
\newcommand\makeinvisible{\let\hide\hide@invisible}
\makeatother
\makehidden

\usepackage{lipsum}
\usepackage{amsfonts}
\usepackage{graphicx}
\usepackage{epstopdf}
\usepackage{algorithmic}
\ifpdf
  \DeclareGraphicsExtensions{.eps,.pdf,.png,.jpg}
\else
  \DeclareGraphicsExtensions{.eps}
\fi

\newsiamremark{remark}{Remark}
\newsiamremark{hypothesis}{Hypothesis}
\crefname{hypothesis}{Hypothesis}{Hypotheses}
\newsiamthm{claim}{Claim}
\newsiamremark{fact}{Fact}
\crefname{fact}{Fact}{Facts}

\headers{Efficient BEM for the Smoluchowski Diffusion Equation}{I.~Labarca-Figueroa, H.~Gimperlein}

\title{Efficient Boundary Elements for the Smoluchowski Diffusion Equation\thanks{
\funding{This research was funded in part by the Austrian Science Fund (FWF) No. 10.55776/P35673.}}}

\author{Ignacio Labarca-Figueroa\textsuperscript{$\dag$,$\ddag$,$\S$}
\and Heiko Gimperlein\thanks{Engineering Mathematics, University of Innsbruck (\email{heiko.gimperlein@uibk.ac.at}).} 
}

\usepackage{amsopn}

\ifpdf
\hypersetup{
  pdftitle={Boundary Integral Formulations for Diffusion under Linear Shear Flow},
  pdfauthor={Ignacio Labarca-Figueroa, et al.}
}
\fi

\begin{document}

\maketitle

\begingroup
\renewcommand{\thefootnote}{\fnsymbol{footnote}}
\footnotetext[3]{Institute for Theoretical Physics, University of Innsbruck.}
\footnotetext[4]{Institute for Mathematical and Computational Engineering, School of Engineering and Faculty of Mathematics, Pontificia Universidad Cat\'olica de Chile, Santiago, Chile (\email{ignacio.labarca@uc.cl}).}
\endgroup

% REQUIRED
\begin{abstract}
The Smoluchowski diffusion equation describes diffusion in the presence of external forces. Studying the mechanical response of soft materials to linear forces, such as shear, results in a boundary value problem involving an Ornstein-Uhlenbeck operator in an exterior domain with non-constant, unbounded coefficients. In this article, we present efficient and highly accurate boundary element methods in the frequency domain, motivated by applications in soft matter physics. Our key contributions concern the accurate assembly of the Galerkin matrix, combining the approximation of the fundamental solution as a Fourier integral with the resolution of near-field singularities. Numerical experiments demonstrate the accuracy and efficiency of the proposed methods and show their relevance for the computation of rheological quantities. 
\end{abstract}

% REQUIRED
%\begin{keywords}
%boundary element method, approximate fundamental solution, Smoluchowski equation, soft matter
%\end{keywords}
% REQUIRED
%\begin{MSCcodes}
%68Q25, 68R10, 68U05
%\end{MSCcodes}

\section{Introduction}\label{sec:intro}

Soft matter is a class of materials which are neither liquid nor solid, but can be easily deformed. It includes a wide range of biological tissues,  suspensions, foams and granular dispersions. A central question concerns their  flow under imposed forces, and in particular shear.

Colloidal suspensions, like biological membranes, sediments or also mayonnaise, consist of typically sub-micrometer-sized particles suspended in a liquid. They offer a model system in which Brownian diffusion, hydrodynamic advection and interparticle forces compete \cite{dhont1996introduction}. While shear stress is proportional to the applied shear rate in Newtonian fluids, soft materials may exhibit a nonlinear response due to these competing effects \cite{brader2010nonlinear}. 
Unfortunately, it is often analytically intractable to study the macroscopic rheology resulting from a given microscopic structure. 
In the simplest such problem, an advection-diffusion equation named after Smoluchowski, 
\begin{equation}\label{eq:smoluchowski_equation_r_intro}
 \partial_t \psi -D_0\Delta_{\nex} \psi  +   \nabla_{\nex} \cdot (\nev\psi)=  0,\end{equation}
 models the evolution of the probability density of a colloidal particle in a suspension under a given flow \cite{brader2010nonlinear}.
 Here, $D_0>0$ is the diffusivity constant, and for a linear shear flow the vector field $\nev$ is given by $\nev(\nex) = \dot{\gamma} x_2 \hat{\nex}_1$% = \Pe B\nex$
 , where $\dot{\gamma} > 0$ denotes the shear rate.
More generally, the shear response of soft materials at low densities is described by the Smoluchowski diffusion equation \eqref{eq:smoluchowski_equation_r_intro} for pairs of particles, to leading order in the occupied volume fraction \cite{bergenholtz2002non,brady1995normal}. A central goal is to compute rheological quantities, i.e.~macroscopic response functions of the suspension to an applied deformation, such as viscosity and normal stress differences \cite{blawzdziewicz1993structure,brader2010nonlinear, brady1995normal}. It is worth noting that \eqref{eq:smoluchowski_equation_r_intro} has numerous further applications as the Fokker-Planck equation of the  Ornstein-Uhlenbeck stochastic process \cite{gardiner}, ranging from financial mathematics to evolutionary biology.

This article investigates efficient and highly accurate numerical methods for solving initial-boundary value problems for the Smoluchowski diffusion equation. Of particular interest is the case of linear shear flow in unbounded exterior domains $\Omega$,  for which previous fundamental studies include  \cite{bergenholtz2002non, brady1995normal, blawzdziewicz1993structure, squires2005simple}. The achieved accuracy of the developed methods will allow in \cite{future} to study new aspects of the rheology.

We first reduce the boundary value problems for the Smoluchowski equation \eqref{eq:smoluchowski_equation_r_intro} from the domain $\Omega$ to the boundary $\partial \Omega$, in the frequency domain. The Galerkin approximation of the resulting boundary integral equations leads to efficient boundary element methods (BEM) for their numerical solution. A central algorithmic contribution of this article is the accurate assembly of the Galerkin matrix, where we handle the numerical approximation of the fundamental solution by windowed truncations of a Fourier integral and careful resolution of the near-field singularities. Numerical experiments investigate the accuracy of evaluating the matrix entries and the convergence of the numerical solutions. Finally, we consider the resulting numerical approximations of rheological quantities obtained as functionals of the boundary element solutions.

Boundary element methods offer several advantages over finite element
methods (FEM) in exterior boundary value problems for \eqref{eq:smoluchowski_equation_r_intro}.
Standard FEM approaches require  the computational domain $\Omega$ to be artificially truncated and radiation conditions to be introduced 
at infinity; additionally, the advection-dominated regime at large P\'eclet numbers
leads to boundary layers and highly anisotropic solutions that require either fine
meshes or stabilization techniques. We also note the unbounded coefficients $\nev(\nex)$ for the shear problem.  By contrast, the BEM reduces the problem to the
boundary $\partial\Omega$ alone, automatically handles the behavior at infinity and unbounded coefficients, and
avoids volumetric discretization errors. 

Classical references for the BEM include the monographs \cite{sauter2010boundary, steinbach2008numerical}, while the underlying analysis in Lipschitz domains is carried out in \cite{costabel1988boundary, mclean2000strongly}. It relies on the availability of fundamental solutions.
For the Smoluchowski equation under linear shear, the time-domain fundamental
solution is known in closed form as a modified anisotropic Gaussian
(see Section~\ref{sec:repformula}). However, for stationary and frequency-domain
problems, one requires its Fourier--Laplace transform
%, which does not admit a closed form. 
. As mentioned above, accurately evaluating this transform numerically---handling the near-field
singularity at small times, the algebraic tail decay, and for positive frequencies the
oscillatory behavior of the integrand---is the central algorithmic contribution of this article. Related strategies for approximating or constructing fundamental
solutions in non-standard settings include the method of
images~\cite{cheng1998method}, efficient representations for Stokes
flow~\cite{broms2025method, gimbutas2015simple}, Ewald summation and quasi-periodic Green functions for periodic
problems~\cite{arens2013analysing, bruno2014rapidly, lindbo2011spectral}, the
windowed Green function method for layered-media
geometries~\cite{bruno2016windowed, bruno2017windowed}, explicit factorizations for
impedance half-space
problems~\cite{lin2024greensI, lin2024greensII}, and low-rank
approximations~\cite{liang2025accelerating}.

Finally, we note that the frequency domain, which we consider, is of particular interest to physicists to understand the long-time
behavior of solutions to \eqref{eq:smoluchowski_equation_r_intro}: low frequencies correspond to the large-time asymptotics of
the pair distribution function. This has been exploited in the context of active
microrheology, where the response
of a suspension to a probe particle pulled by a constant
force, $\nev(\nex) = \mu \neF$, is studied \cite{franosch2018time,squires2005simple}. This simpler problem, however, can be reduced to a standard heat equation by a change of variables, and appropriate numerical methods are well studied.

\section{Formulation of shear problems for soft materials}\label{sec:physmot}

\begin{figure}
    \centering
    \includegraphics[width=0.90\linewidth]{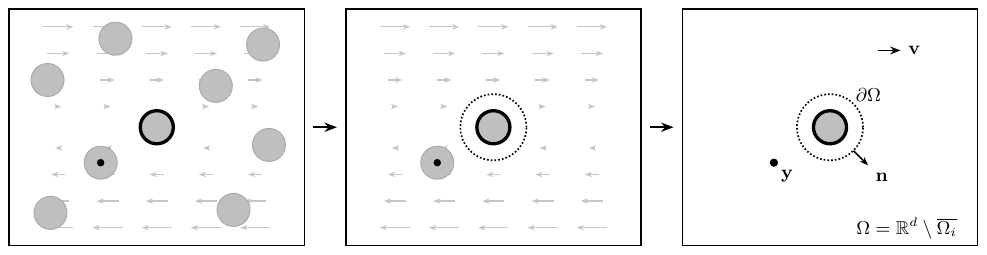}
    \caption{Left: dilute colloidal suspension with several particles. Center: two-particle setting under shear flow. Right: geometry and notation for \eqref{eq:smoluchowski_equation_r}--\eqref{eq:smoluchowski_bc}. } 
    \label{fig:colloid}
\end{figure}

A colloidal suspension consists of mesoscopic particles dispersed in a solvent. Under shear flow, as shown in Figure~\ref{fig:colloid}, left and center, Brownian diffusion and hydrodynamic advection compete to determine the spatial arrangement of the particles. This results in a non-equilibrium microstructure and non-Newtonian behavior, meaning the fluid's viscosity and stress response depend nonlinearly on the applied shear rate \cite{brader2010nonlinear}.

In the dilute limit, the particle volume fraction is small and only pairwise interactions between particles need to be considered. The suspension microstructure is captured by the pair distribution function $\psi(\nex,t\,|\,\ney)$, which is the probability density of finding two colloidal particles at relative position $\nex$ at time $t$, given that they started at $\ney$ (see Figure~\ref{fig:colloid}, right). To leading order in the volume fraction, the time evolution of $\psi$ is governed by the Smoluchowski diffusion equation
\begin{equation*}
  \partial_t \psi - D_0\Delta\psi + \nabla\cdot(\nev\psi) = 0 \quad \text{in } \mathbb{R}^3\setminus\overline{\Omega_i},\ t > 0,
\end{equation*}
where $\Omega_i$ is the excluded volume; for spherical particles of radius $\sigma$, $\Omega_i = \{\nex : |\nex| < 2\sigma\}$ is a sphere of radius $2\sigma$, the contact distance of two particles.  The associated Neumann boundary condition,
\begin{equation*}
  \gammas \psi \coloneqq \mathbf{n}\cdot(\nev\psi - D_0\nabla\psi) = 0 \quad \text{on } \partial\Omega_i,
\end{equation*}
expresses a no-flux condition corresponding to hard-sphere repulsion: no probability flux passes through the contact surface. For linear shear flow $\nev(\nex) = \dot{\gamma}\,x_2\,\hat{\nex}_1$, the P\'eclet number $\Pe = \dot{\gamma}\sigma^2/D_0$ measures the relative importance of advection and Brownian diffusion. At equilibrium ($\Pe = 0$), the pair distribution is isotropic with $\psi \to 1$ as $|\nex|\to\infty$; for $\Pe > 0$, shear breaks this isotropy and the microstructure becomes anisotropic.\\

Rheology studies the flow and deformation properties of complex fluids. The macroscopic transport coefficients of interest (the shear viscosity, the normal stress differences, and the frequency-dependent viscosity) characterize non-Newtonian phenomena such as shear thinning, shear thickening, and viscoelasticity \cite{bergenholtz2002non,brader2010nonlinear}. These quantities are all functionals of $\psi$, see \cite{blawzdziewicz1993structure,brader2010nonlinear,brady1995normal,squires2005simple} for the stress tensor and \cite{brader2010nonlinear} for the stress autocorrelation function. In particular, the time-dependent (averaged) stress tensor $\langle \sigma^{\alpha\beta} \rangle_t$
is given by~\cite{brader2010nonlinear}
\begin{equation}\label{eq:stress-tensor}
    \langle \sigma^{\alpha\beta} \rangle_t = -nk_BT\delta^{\alpha\beta}
    + \frac{5}{2}\eta\varphi(\kappa^{\alpha\beta} + \kappa^{\beta\alpha})
    - \frac{n\sigma^3}{2}nk_BT \int_{\partial\Omega} 
    \frac{x^\alpha x^\beta}{|\nex|^2}\, \Psi(\nex, t)\dsx,
\end{equation}
where $$\Psi(\nex, t) = \int_{\Omega} \psi(\nex, t \mid \ney)\, \dy$$ is the pair
correlation function.
The three contributions correspond, respectively, to the ideal gas pressure, the Einstein
correction from the applied shear %$\nev(\nex) = \Pe\, x_2\, \hat{\nex}_1$ 
with velocity
gradient tensor $\kappa^{\alpha\beta} = \Pe\, \delta_{\alpha 1}\delta_{\beta 2}$, and the
hard-sphere repulsion at contact. Here $\eta$ is the solvent viscosity and
$\varphi = \pi n \sigma^3 / 6$ is the packing fraction of hard spheres with number
density $n = N/V$ and diameter $\sigma$.

\section{Smoluchowski diffusion equation and boundary integral formulation}\label{sec:repformula}
%\section{Time-Domain Representation formula, boundary integral formulation} 

\subsection{Time- and frequency-dependent boundary value problems}

Let  $\Omega_i$ be a bounded Lipschitz domain in $\mathbb{R}^d$ with boundary $ \Gamma \coloneqq \partial \Omega$. As motivated in Section \ref{sec:physmot}, we consider the  initial-boundary value problem for the Smoluchowski diffusion equation \eqref{eq:smoluchowski_equation_r_intro} with diffusion coefficient $D_0=1$ in the exterior  $\Omega \coloneqq \mathbb{R}^d\setminus \overline{\Omega_i}$: Find the density $\psi = \psi(\nex, t \, |\, \ney)$ such that
\begin{align}\label{eq:smoluchowski_equation_r}
&\text{PDE:} &
 \partial_t \psi -\Delta \psi  +   \nabla \cdot (\nev\psi)&=  0  && \text{for } \nex \in \Omega,\ t>0,\\ \label{eq:initialcond}
&\text{Initial condition:} &  \psi(\nex, t = 0\, |\, \ney) &=\delta(\nex-\ney)  && \text{for } \nex \in \Omega,\\ \label{eq:smoluchowski_bc}
&\text{Boundary condition:} &\mathbf{n} \cdot \left( \nev \psi -  \nabla \psi \right) &= 0 && \text{for } \nex \in \partial \Omega,\ t>0,
\end{align}
for each $\ney \in \Omega$. Here, we assume that $\psi$ is bounded as $|\nex| \to \infty$. We denote the unit normal pointing outside of $\Omega_i$ by $\mathbf{n}$ and consider velocities $\nev$ that correspond to an incompressible flow, i.e. $\nabla \cdot \nev = 0$. As discussed in Section \ref{sec:intro}, linear shear flows are of particular interest, corresponding to $\nev(\nex) = \Pe x_2 \hat{\nex}_1$% = \Pe B\nex$
, with P\'eclet number $\Pe > 0$.

%The velocity $\nev$ is assumed to be a (possibly non-constant) incompressible flow, i.e. $\nabla \cdot \nev = 0$. 
%We consider two examples: 
%\begin{enumerate}[label=(\roman*)]
%    \item $\nev(\nex) = \mu \neF$, where $\mu > 0$ is the mobility and $\neF(\nex) \equiv \neF$ is a constant vector field.
%    \item $\nev(\nex) = \Pe x_2 \hat{\nex}_1 = \Pe B\nex$, where $\Pe > 0$ is the P\'eclet number of the problem. \\ We write
%    $$ B = \begin{pmatrix}
%        0 & 1 & 0\\ 0 & 0 &0 \\ 0& 0 & 0
%    \end{pmatrix},$$
%    then $\nev$ corresponds to a linear shear-flow.
%\end{enumerate}

For the numerical solution, it is convenient to translate this problem into a problem with zero initial conditions, but inhomogeneous boundary conditions, using the free-space solution $\psi_0$ in $\mathbb{R}^d$:
\begin{align}\label{eq:smoluchowski_equation_free}
\begin{aligned}
&\text{PDE:} & 
 \partial_t \psi_0 -\Delta \psi_0 +\nabla \cdot (\nev\psi_0)&=0    && \text{for } \nex \in \mathbb{R}^d,\ t>0,\\
&\text{Initial condition:} &  \psi_0(\nex, t = 0\, |\, \ney) &=\delta(\nex-\ney)  &&\text{for } \nex \in \mathbb{R}^d.
\end{aligned}
\end{align}
In fact, from Section \ref{sec:physmot} we recall that the quantities of physical interest can be computed from the expectation value $\Psi$ in $\ney$, $$\Psi(\nex, t) = \displaystyle\int\limits_{\Omega} (\psi(\nex, t \, | \, \ney) - \psi_0(\nex, t \, | \, \ney))\dy = \displaystyle\int\limits_{\Omega} \psi(\nex, t \, | \, \ney)\dy -1 $$ satisfies the initial-boundary value problem
\begin{align}\label{eq:smoluchowski_equation_Psi}
&\text{PDE:} & 
 \partial_t \Psi -\Delta \Psi + \nabla \cdot(\nev \Psi)  &= 0&& \text{in } \Omega,\ t>0, \\ \label{eq:initialcond_Psi}
&\text{Initial condition:} &  \Psi(\nex, t = 0) &=0  && \text{in }  \Omega,\\
&\text{Boundary condition:} & \gammas \Psi &= g_s    && \text{on }  {\partial \Omega},\ t>0. \label{eq:smoluchowski_equation_Psibc}
\end{align}
Here, the  boundary data are given by $g_s(\nex) = -\mathbf{n} \cdot \nev(\nex)$, and $\Psi$ tends to $0$ as $|\nex| \to \infty$. The analysis of the model problems of constant force or constant shear is particularly well understood, as special cases of the heat and the Ornstein-Uhlenbeck operators.\\

As a linear parabolic problem, we may solve \eqref{eq:smoluchowski_equation_Psi}-\eqref{eq:smoluchowski_equation_Psibc} using the Fourier-Laplace transform or, equivalently, using a time-harmonic ansatz with frequency $\omega\geq0$,
\begin{equation}\label{eq:time-harmonic}
\Psi(\nex, t) = \hat{\Psi}_\omega(\nex) \exp(i\omega t).
\end{equation}
The frequency-domain solution $\hat{\Psi}_\omega$ satisfies the \emph{exterior} boundary value problem
\begin{align}\label{eq:smoluchowski_equation_Psi-freq}
\begin{aligned}
&\text{PDE} & 
 i\omega \hat{\Psi}_\omega -\Delta \hat{\Psi}_\omega +  \nabla \cdot(\nev \hat{\Psi}_\omega)&= 0 && \text{for } \nex \in \Omega, \\
&\text{Boundary condition } & \gamma_{\star}\hat{\Psi}_\omega &= g_{\star}   && \text{for } \nex \in {\partial \Omega}. 
\end{aligned}
\end{align}
Here, $\star  =\mathrm{s}$ for the Fourier-Laplace transform of the physical (Neumann) problem \eqref{eq:smoluchowski_equation_Psi}-\eqref{eq:smoluchowski_equation_Psibc}. Below, however, we also consider the Dirichlet problem with $\star  =\mathrm{D}$, as well as the corresponding \emph{interior} Dirichlet and Neumann problems, where $\nex \in \Omega_i$. 
\subsection{Fundamental solutions in time and frequency domains}
Boundary integral formulations for the solution of inhomogeneous boundary value problems rely on the knowledge of a fundamental solution for the PDE. In the case of the time-dependent model problem \eqref{eq:smoluchowski_equation_r} for constant shear  in $\mathbb{R}^3$, the 
%\begin{equation}
%     \partial_t \Psi -\Delta \Psi + \nabla \cdot(\nev \Psi)  = 0 \qquad \text{in } \mathbb{R}^3,\ t>0,
%\end{equation}
%it is possible to write the 
fundamental solution $G$ is given by a modified Gaussian function \cite{blawzdziewicz1993structure,elrick1962source}: When $\Pe = 1$, we have
\begin{equation}\label{eq:fundamental-time1}
\begin{aligned}
    G(\nex, \ney, t) &= \tfrac{1}{(4\pi t)^{3/2} \sqrt{1+t^2 / 12}} \exp\left(-\frac{p(\nex, \ney, t)^2}{4t}\right),\\
    p(\nex, \ney, t)^2 &=  \tfrac{(x_1 - y_1 - t(x_2 + y_2)/2)^2}{1+t^2 / 12} +  (x_2 -y_2)^2 + (x_3 - y_3)^2.
\end{aligned}
\end{equation}
The fundamental solution for $\Pe\neq 1$ is obtained by rescaling in space and time, % the spatial and time variables,
\begin{equation*}
    \nex = \Pe^{1/2} \nex^{\star}, \quad t = \Pe t^{\star}, \quad \psi = \Pe^{-3/2}\psi^{\star}. 
\end{equation*}
In particular, in  Section \ref{sec:repformula-frequency} we use the following formula in real space variables, but scaled time variables:
\begin{align}\label{eq:fundamental-solution-shear}
G(\nex^{\star}, \ney^{\star}, \tau) &= \tfrac{1}{(4\pi \tau)^{3/2} \sqrt{1+\tau^2 / 12}}\exp\left(-\Pe \tfrac{p(\nex^{\star}, \ney^{\star}, t)^2}{4\tau} \right).
\end{align}
Using the Fourier-Laplace transform, from \eqref{eq:fundamental-time1}, we also obtain an integral formula for the fundamental solution $\hat{G}_{\omega}$ of the PDE \eqref{eq:smoluchowski_equation_Psi-freq} in $\mathbb{R}^3$ in the frequency domain:
\begin{equation}\label{eq:fundamental-solution-shear-freq}
    \hat{G}_{\omega}(\nex, \ney) = \displaystyle \int_0^{\infty} G(\nex, \ney, t) \exp(-i\omega t)dt. 
\end{equation}

\subsection{Representation formula, boundary integral formulation} \label{sec:repformula-frequency}

The fundamental solution \eqref{eq:fundamental-solution-shear-freq} allows us to define the standard layer potentials and boundary integral operators in the frequency domain. 
A solution of the differential equation $$i\omega \hat{\Psi}_\omega -\Delta \hat{\Psi}_\omega +  \nabla \cdot(\nev \hat{\Psi}_\omega)= 0$$ in $\mathbb{R}^3 \setminus \Gamma$, which decays as $|\nex|\to \infty$, satisfies the representation formula
\begin{equation}\label{eq:repformula-freq}\hat{\Psi}_\omega(\nex) =  \hat{\mathrm{D}}_{\omega}\jump{\gammad\hat{\Psi}_\omega}(\nex) - \hat{\mathrm{S}}_{\omega}\jump{\gammas\hat{\Psi}_\omega}(\nex).\end{equation}
Here, $\jump{\cdot}$ denotes the jump across the boundary $\Gamma$, 
\begin{equation}\label{eq:jump-op}
    \jump{\gamma_{\star} \hat{\Psi}_\omega}  \coloneqq \gamma_{\star}^+ \hat{\Psi}_\omega - \gamma_{\star}^- \hat{\Psi}_\omega, \qquad \text{ for }\star \in \{\mathrm{s}, \mathrm{D}\},
\end{equation}
where $\gamma_{\star}^+$ (resp. $\gamma_{\star}^-$) denotes the trace taken from the exterior (resp. interior) of $\Omega_i$.\\
For smooth densities $\nu$ and $\phi$ the single and double layer potentials are given by
\begin{align}\label{eq:layer-potentials}
    \hat{\mathrm{S}}_{\omega}\nu (\nex) &= \int_{\partial \Omega} \hat{G}_{\omega}(\nex,\ney) \nu(\ney) \dsy,  \\
    \hat{\mathrm{D}}_{\omega}\phi (\nex) &= \int_{\partial \Omega} \dfrac{\partial \hat{G}_{\omega}}{\partial \mathbf{n}_{\ney}}(\nex,\ney) \phi(\ney) \dsy.
\end{align} 
We introduce the boundary integral operators (BIOs)
\begin{equation}\label{eq:BIOs}
\begin{aligned}
&\hat{\mathrm{V}}_{\omega} = \average{\gammad \hat{\mathrm{S}}_{\omega}}, \qquad && \hat{\mathrm{K}}_{\omega} = \average{\gammad\hat{\mathrm{D}}_{\omega}}, \\ 
 &\W_{\omega} = \average{\gammas\hat{\mathrm{D}}_{\omega}}, \qquad &&
\hat{\widetilde{\mathrm{K}}}^{'}_{\omega} = \average{\gammas\hat{\mathrm{S}}_{\omega}},
\end{aligned}    
\end{equation}
where $\average{ \cdot }$ denotes the average across the boundary $\Gamma$:
\begin{equation}
    \average{\gamma_{\star}\Phi} \coloneqq \frac{1}{2}\left( \gamma_{\star}^+ \Phi + \gamma_{\star}^- \Phi\right), \qquad \text{ for }\star \in \{\mathrm{s}, \mathrm{D}\}.
\end{equation}
The following jump relations for the layer potentials are classical \cite[Lemma~4.1]{costabel1988boundary}:
\begin{equation}\label{eq:jump-condition}
\begin{aligned}
\jump{\gammad\mathrm{S}_{\omega}\nu} &= 0,&\quad\jump{\gammas\mathrm{S}_{\omega}\nu)} &= -\nu, \\
\jump{\gammad\mathrm{D}_{\omega}\phi} &= \phi,&\quad\jump{\gammas\mathrm{D}_{\omega}\phi)} &= 0.    
\end{aligned}
\end{equation}

Using the representation formula \eqref{eq:repformula-freq}, the definition of BIOs \eqref{eq:BIOs}, and the jump relations \eqref{eq:jump-condition}, we now reformulate the Dirichlet and Neumann boundary value problems \eqref{eq:smoluchowski_equation_Psi-freq} as equivalent boundary integral equations (BIEs). \\

To be specific, for the \emph{exterior} problems in $\Omega$, $\hat{\Psi}_{\omega} \equiv 0$ in ${\Omega_i}$, and the representation formula \eqref{eq:repformula-freq} simplifies to
\begin{equation}\label{eq:repformula-exterior}
    \hat{\Psi}_{\omega}(\nex) = \hat{\mathrm{D}}_{\omega}(\gamma^+_{\mathrm{D}}\hat{\Psi}_{\omega})(\nex) - \hat{\mathrm{S}}_{\omega}(\gamma^+_{\mathrm{s}}\hat{\Psi}_{\omega})(\nex).\end{equation}
The jump relations \eqref{eq:jump-condition}  then lead to the following BIEs:
\begin{proposition}
Let $\omega \geq 0$ and $\star \in \{\mathrm{s}, \mathrm{D}\}$. Then the \emph{exterior} boundary value problem \eqref{eq:smoluchowski_equation_Psi-freq} is equivalent to 
\begin{align}\label{eq:BIEs_ext}
\begin{aligned}
    &\text{SL Direct: } & \hat{\mathrm{V}}_{\omega} \psi &= \left( -\tfrac{1}{2}\mathrm{I} + \hat{\mathrm{K}}_{\omega} \right)g_D & \text{($\star = \mathrm{D}$, Dirichlet Problem)},\\
    &\text{DL Direct: } & \left( -\tfrac{1}{2}\mathrm{I} + \hat{\mathrm{K}}_{\omega} \right)\varphi &= \hat{\mathrm{V}}_{\omega} g_{\mathrm{s}} & \text{($\star = \mathrm{s}$, Neumann Problem)},
\end{aligned}
\end{align}
with $\psi = \gamma^+_{\mathrm{s}}\hat{\Psi}_{\omega}$, $\varphi =  \gamma^+_{\mathrm{D}}\hat{\Psi}_{\omega}$. 
\end{proposition}
The precise function spaces, in which these integral equations are wellposed, will be reviewed from \cite{costabel1988boundary} in Section \ref{sec:BEM}.

Similarly, for the \emph{interior} problems in $\Omega_i$, $\hat{\Psi}_\omega \equiv 0$ in $\Omega$, and the representation formula \eqref{eq:repformula-freq} simplifies to
\begin{equation}\label{eq:repformula-interior}
    \hat{\Psi}_\omega(\nex) = \hat{\mathrm{S}}_{\omega}(\gamma^-_{\mathrm{s}}\hat{\Psi}_\omega)(\nex) - \hat{\mathrm{D}}_{\omega}(\gamma^-_{\mathrm{D}}\hat{\Psi}_\omega)(\nex).\end{equation}
The jump relations \eqref{eq:jump-condition}  then lead to:
\begin{proposition}
Let $\omega \geq 0$ and $\star \in \{\mathrm{s}, \mathrm{D}\}$. Then the \emph{interior} boundary value problem \eqref{eq:smoluchowski_equation_Psi-freq} is equivalent to 
\begin{align}\label{eq:BIEs}
\begin{aligned}
    &\text{SL Direct: } & \hat{\mathrm{V}}_{\omega} \psi &= \left( \tfrac{1}{2}\mathrm{I} + \hat{\mathrm{K}}_{\omega} \right)g_{\mathrm{D}} & \text{($\star = \mathrm{D}$, Dirichlet Problem)},\\
    &\text{DL Direct: } & \left( \tfrac{1}{2}\mathrm{I} + \hat{\mathrm{K}}_{\omega} \right)\varphi &= \hat{\mathrm{V}}_{\omega} g_{\mathrm{s}} & \text{($\star = \mathrm{s}$, Neumann Problem)}.
\end{aligned}
\end{align}
with $\psi = \gamma^-_{\mathrm{s}}\hat{\Psi}_{\omega}$, $\varphi = \gamma^-_{\mathrm{D}}\hat{\Psi}_{\omega}$. 
\end{proposition}

\remark{As an alternative to the representation formula \eqref{eq:repformula-freq}, one may consider a layer potential ansatz:
\begin{equation}\label{eq:rep-formula-indirect}
\begin{aligned}
    &&\hat{\Psi}_\omega(\nex) = (\hat{\mathrm{S}}_{\omega}\lambda)(\nex), \qquad \nex \in \mathbb{R}^3 \setminus \Gamma,\\
    \text{or} &&\hat{\Psi}_\omega(\nex) = (\hat{\mathrm{D}}_{\omega}\phi)(\nex), \qquad \nex \in \mathbb{R}^3 \setminus \Gamma.
\end{aligned}
\end{equation}
The \emph{exterior} Dirichlet problem \eqref{eq:smoluchowski_equation_Psi-freq} with $\star = \mathrm{D}$ is then equivalent to
\begin{align}\label{eq:BIEs-indirect_ext}
\begin{aligned}
    &\text{SL Indirect: } & \hat{\mathrm{V}}_{\omega} \lambda &= g_{\mathrm{D}} & \text{($\star = \mathrm{D}$, Dirichlet Problem)},\\
    &\text{DL Indirect: } & \left( \tfrac{1}{2}\mathrm{I} + \hat{\mathrm{K}}_{\omega} \right)\phi &= g_{\mathrm{D}} & \text{($\star = \mathrm{D}$, Dirichlet Problem)},
\end{aligned}
\end{align}
while the \emph{interior} Dirichlet problem is equivalent to
\begin{align}\label{eq:BIEs-indirect}
\begin{aligned}
    &\text{SL Indirect: } & \hat{\mathrm{V}}_{\omega} \lambda &= g_{\mathrm{D}} & \text{($\star = \mathrm{D}$, Dirichlet Problem)},\\
    &\text{DL Indirect: } & \left( -\tfrac{1}{2}\mathrm{I} + \hat{\mathrm{K}}_{\omega} \right)\phi &= g_{\mathrm{D}} & \text{($\star = \mathrm{D}$, Dirichlet Problem)}.
\end{aligned}
\end{align}
}

\section{Numerical approximation of frequency-domain fundamental solutions}\label{sec:approx}
In this section we describe the numerical computation of the fundamental solution its use in the computation of the boundary integral operators. From \eqref{eq:fundamental-time1} we can point out some properties about the (time-domain) fundamental solution that will impact the approximation of its  counterpart in the frequency domain. 

First, the singular behavior at $\tau=0$ along $\nex = \ney$: there, the fundamental solution exhibits the same singularities as the fundamental solution for the heat equation:
\begin{equation}\label{eq:GH} 
    G_{\text{H}}(\nex, \ney,\tau) \coloneqq \dfrac{1}{(4\pi \tau)^{3/2}}\exp\left(-\tfrac{\Pe}{4\tau}|\nex - \ney|^2 \right).
\end{equation}
Second, the decay at $\tau=\infty$: While for $\nex \neq \ney$, the fundamental solution decays exponentially in $\tau$, the decay becomes algebraic $\mathcal{O}(\tau^{-5/2})$ near $\nex = \ney$.\\

This singular behavior needs to be taken into account in the numerical approximation of the frequency-domain fundamental solution \eqref{eq:fundamental-solution-shear-freq} 
\begin{equation*}
    \hat{G}_{\omega}(\nex, \ney) = \displaystyle \int_0^{\infty} G(\nex, \ney, t) \exp(-i\omega t)dt. 
\end{equation*}
from the known time-domain fundamental solution \eqref{eq:fundamental-solution-shear}. 
For an accurate numerical quadrature, we split the time integral into two terms. We fix $\tau_0 > 0$ and write
\begin{align}\label{eq:split}
    \hat{G}_{\omega}(\nex, \ney) &= \underbrace{\displaystyle \int_0^{\tau_0} G(\nex, \ney, t) \exp(-i\omega t)dt}_{\Ising(\nex, \ney)} 
    +\underbrace{\displaystyle \int_{\tau_0}^{\infty} G(\nex, \ney, t) \exp(-i\omega t)dt}_{\Ilong(\nex, \ney)}. 
\end{align}
For any fixed $\nex, \ \ney \in \mathbb{R}^d$ the fundamental solution for the shear problem in the time domain decays (at least) algebraically fast in time. This algebraic decay ensures the convergence of $\Ilong$ and suggests adequate quadrature rules for the approximation. The approaches differ for the cases of positive ($\omega>0$) and zero ($\omega=0$) frequencies.

\remark{The numerical quadrature will be applied with a rescaled time variable $\frac{t}{\Pe}$. In applications the corresponding rescaled frequency $\omega \Pe$ is typically boun\-ded or even small. } %Recall that we use rescaled time-variables. Therefore, the real frequency corresponds to $\omega \Pe$.  For example, if $\Pe = \tfrac{1}{10},$ and we use $\omega = \tfrac{1}{2}$, then the real frequency is $\tfrac{1}{20}$. }
\subsection{Stationary problem ($\omega = 0$)}
\subsubsection{Evaluating $\Ising$}
When $\omega  = 0$, we rely on the following result for $\hat{G}_{\text{Lap}}(\nex, \ney)$, a multiple of the fundamental solution of the Laplacian 
\begin{equation}
    \hat{G}_{\text{Lap}}(\nex, \ney) = \dfrac{1}{4\pi{\Pe^{1/2}}|\nex - \ney|} = \int_0^{\infty} \underbrace{\dfrac{1}{(4\pi \tau)^{3/2}}\exp\left(-\tfrac{\Pe}{4\tau}|\nex - \ney|^2 \right)}_{G_{\text{H}}(\nex, \ney,\tau)}d\tau.
\end{equation}
More generally, for a time integral over $(0,\tau_0)$
\begin{equation}
    \dfrac{\erfc\left(\tfrac{{\Pe^{1/2}}}{2\sqrt{\tau_0}}|\nex - \ney|\right)}{4\pi\Pe^{1/2}|\nex - \ney|} = \int_0^{\tau_0} \dfrac{1}{(4\pi \tau)^{3/2}}\exp\left(-\tfrac{\Pe}{4\tau}|\nex - \ney|^2 \right)d\tau,
\end{equation}
where $\erfc$ is the complementary error function
$    \erfc(x) = 1- \dfrac{2}{\sqrt{\pi}}\int_0^x \exp(-x^2)dx$. We define $\Gdiff$ as follows:
\begin{equation}\label{eq:Gdiff}
    \Gdiff(\nex, \ney,\tau) \coloneqq G(\nex, \ney,\tau) - G_{\text{H}}(\nex, \ney,\tau) \exp\left(  \tfrac{\Pe}{4}(x-x')(y+y')\right).
\end{equation}
Algebraic manipulations lead to the following representation of $\Gdiff$,
\begin{equation}\label{eq:GHf}
    \Gdiff(\nex, \ney,\tau) = G_{\text{H}}(\nex, \ney,\tau) f(\nex, \ney, \tau),
\end{equation}
where
\begin{equation}\label{eq:continuous}
    f(\nex, \ney, \tau) \coloneqq \dfrac{\exp\left(\tau\Pe\tfrac{(x-x')^2/12 -(y+y')^2/4 - \tau(x-x')(y+y')/12}{4(1+\tau^2/12)}\right)}{\sqrt{1+\tau^2/12}} - 1.
\end{equation}
The function $f$ is $\mathcal{O}(\tau)$ for small values of $\tau$, uniformly in $\nex, \ney$.
For $|\nex - \ney| > 0$, $G_{\text{H}} \rightarrow{0}$ exponentially fast as $\tau\rightarrow 0$, and therefore so does $\Gdiff$. \\
For $|\nex - \ney| = 0$, $G_{\text{H}}$ behaves as $\mathcal{O}(\tau^{-3/2})$ for small times, i.~e.,~it is not integrable, so that with \eqref{eq:GHf} we see that $\Gdiff$ behaves as $\mathcal{O}(\tau^{-1/2})$ and is integrable.\\

We conclude that 
\begin{equation}\begin{aligned}
    \Isingz(\nex, \ney) = \int_0^{\tau_0} G(\nex, \ney, t) dt &=   \underbrace{\dfrac{\erfc\left(\tfrac{\Pe^{1/2}}{2\sqrt{\tau_0}}|\nex - \ney|\right)}{4\pi\Pe^{1/2}|\nex - \ney|} \exp\left(  \tfrac{\Pe}{4}(x-x')(y+y')\right)}_{\hat{G}^{\text{sing}}_0(\nex, \ney)} \\ &\qquad + \underbrace{\int_0^{\tau_0} \Gdiff(\nex, \ney, \tau) d\tau}_{\mathcal{I}^{\text{diff}}_0(\nex, \ney)}.\end{aligned}
\end{equation}
The integral $\mathcal{I}_0^{\text{diff}}$ can be efficiently approximated by a Gauss quadrature rule. 

\subsubsection{Evaluating $\Ilongz$}\label{sec:Ilong}
Now we focus on the approximation of the integral
\begin{equation*}
\Ilongz(\nex, \ney) = \int_{\tau_0}^{\infty} G(\nex, \ney, \tau)d\tau.
\end{equation*}
Note that for long times, the fundamental solution for the shear problem behaves as $\mathcal{O}(\tau^{-5/2})$. The change of variables
\begin{equation}\label{eq:change-variables}
    u = \tau^{-3/2}, \qquad du = -\dfrac{3}{2}\tau^{-5/2}d\tau
\end{equation}
leads to a continuous integrand over the finite interval $[0, \tau_0^{-3/2}]$. Gauss quadrature can be used to efficiently approximate $\Ilongz$.

\subsection{Frequency-domain problem ($\omega > 0$)}\label{sec:freq-approx}
For positive $\omega$, additional care is required for the integration of an oscillatory integrand. We focus on the modifications compared to $\omega=0$.
\subsubsection{Evaluating $\Ising$}
When $\omega  > 0$, we denote the positive root of the polynomial $x^4 + \omega^2$ by $\kappa$. We rely on the following result 
for the truncated integral
\begin{equation}
    \hat{G}^{\text{sing}}_{\omega}(\nex, \ney) \coloneqq \int_0^{\tau_0} \dfrac{\exp(-i\omega \tau)}{(4\pi \tau)^{3/2}}\exp\left(-\tfrac{\Pe}{4\tau}|\nex - \ney|^2 \right)d\tau,
\end{equation}
given by
\begin{align}\label{eq:Gsing_w}
    \hat{G}^{\text{sing}}_{\omega}(\nex, \ney) &= \dfrac{\exp(-\kappa\Pe^{1/2}|\nex - \ney|)}{8\pi\Pe^{1/2}|\nex - \ney|}\left(\erf\left(\kappa\sqrt{\tau_0} - \tfrac{\Pe^{1/2}}{2\sqrt{\tau_0}}|\nex - \ney|\right) \right.\\
    &\qquad + \left.\exp(2\kappa\Pe^{1/2}|\nex - \ney|) \erfc\left(\kappa\sqrt{\tau_0} + \tfrac{\Pe^{1/2}}{2\sqrt{\tau_0}}|\nex - \ney|\right) +1\right).
\end{align}
From \eqref{eq:Gdiff}, and \eqref{eq:Gsing_w} it follows the decomposition:
% \begin{equation}
%     G^{\text{diff}}(\nex, \ney,\tau) \coloneqq G(\nex, \ney,\tau) - G_{\text{H}}(\nex, \ney,\tau) \exp\left(  \tfrac{\Pe}{4}(x-x')(y+y')\right),
% \end{equation}
% and conclude the decomposition 
\begin{equation} 
    \Ising =  \hat{G}^{\text{sing}}_{\omega}(\nex, \ney) + \underbrace{\int_0^{\tau_0} \exp(-i\omega\tau) \, \Gdiff(\nex, \ney, \tau) d\tau}_{\Idiff}.
\end{equation}
The integral $\Idiff$ is efficiently approximated by a composite Gauss quadrature rule. 

 \subsubsection{Evaluating $\Ilong$}\label{sec:Ilong-freq}  For $\omega >0$ it is no longer possible to proceed as in Section \ref{sec:Ilong}, because of the oscillatory behavior of the integrand. We want to approximate the integral
\begin{equation}\label{eq:Ilong_freq}
\Ilong(\nex, \ney) = \int_{\tau_0}^{\infty} \exp(-i\omega\tau)\, G(\nex, \ney, \tau)d\tau.
\end{equation}
Windowing schemes provide accurate approximations, and we follow \cite{bruno2016windowed, bruno2017windowed}. They are based on a slow-rise windowing function (see \cite[Section~3]{bruno2016windowed}),
\begin{equation}\label{eq:window}
    \eta(t, \tau_0, \tau_1) = \left\{ 
    \begin{array}{ll}
        1 & \text{ for } t < \tau_0 \\
        \exp(\exp(-\tfrac{1}{x})/(x-1)) & \text{ for } \tau_0 < t < \tau_1,\ x = \dfrac{t - \tau_0}{\tau_1 - \tau_0},\\ 
        0 & \text{ for } t > \tau_1
    \end{array}
    \right.
\end{equation}
that introduces a smooth decay to zero at a given $\tau_1$. To approximate an integral
\begin{equation}\label{eq:integral-long}
    I = \int_1^{\infty} \dfrac{\exp(-i\omega t)}{t^{3/2}}dt,
\end{equation}
we hence do not truncate \eqref{eq:integral-long} as
\begin{equation}\label{eq:integral-long-truncated}
    I_{L} = \int_1^{L} \dfrac{\exp(-i\omega t)}{t^{3/2}}dt, \quad L > 1,
\end{equation}
but rather truncate smoothly using the smooth window function $\eta$,
\begin{equation}\label{eq:Iwindow}
    I_{L,\eta} = \int_1^{L} \dfrac{\exp(-i\omega t)}{t^{3/2}} \eta(t, cL, L)dt, \quad L > 1, c \in (0, 1).
\end{equation}
As a main advantage of the smooth truncation, the error $I-I_{L,\eta}$ tends to $0$ faster than algebraically as $L \to \infty$. This approach is used for both \eqref{eq:Ilong_freq} and its derivatives.

\subsection{Composite Gauss-Legendre quadrature}

After the singularity subtraction, the integrands obtained are only of limited regularity. This leads to limited convergence rates for standard Gauss quadrature rules. In order to recover high-order convergence rates, a composite Gauss quadrature rule is used, as is standard for $hp$--quadrature rules \cite{babuvska1994p, schwab1998hpfem, zhang2021hpquad}.

To be specific, given a function $v : [0, 1]\rightarrow \mathbb{R}$ with limited regularity at $t = 0$, we split the interval $[0,1]$ into $N+1$ subintervals of the form:
\begin{equation}
\begin{aligned}
    [t_{\ell - 1}, t_{\ell}] \subseteq [0, 1], \quad t_{\ell} \coloneq \sigma^{N - \ell} \quad \text{ for }\ell = 0, \ldots, N,\qquad  t_{-1} \coloneq 0,
\end{aligned}
\end{equation}
where $\sigma \in (0, 1)$ is a given parameter. Then, for $\nu\in \mathbb{N}$, we use $(\ell +1)\nu$ Gauss quadrature points for the interval $[t_{\ell-1}, t_{\ell}]$.

\section{Boundary element methods}\label{sec:galerkin}
We discretize the boundary integral formulations using a Galerkin approach. To do so, we consider a boundary mesh $\mathcal{T}_h$, with $h>0$ the mesh size parameter. Finite dimensional subspaces
\begin{align}\label{eq:bem_spaces}
    \begin{aligned}
&P^0_h \subset H^{-1/2}(\Gamma) &&\triangleq \text{ Piecewise constant functions on }\mathcal{T}_h,\\ 
&P^1_h \subset H^{1/2}(\Gamma)  &&\triangleq \text{ Piecewise linear continuous functions on }\mathcal{T}_h,
    \end{aligned}
\end{align}
associated to $\mathcal{T}_h$ are used for each of the relevant trace spaces on $\Gamma$. 

\subsection{Galerkin discretization}\label{sec:BEM}
The Galerkin approximation of the BIEs \eqref{eq:BIEs_ext} to \eqref{eq:BIEs-indirect} is based on their weak formulations. For the single-layer equation, we have:

Given $g_{\mathrm{D}}\in H^{1/2}(\Gamma)$, find $\psi\in H^{-1/2}(\Gamma)$ such that
\begin{equation}\label{eq:SL-weak}
    \begin{aligned}
    (\text{SL Indirect})\qquad &&\langle \hat{\mathrm{V}}_\omega\psi, \psi^* \rangle = \langle g_{\mathrm{D}}, \psi^*\rangle, 
\end{aligned}
\end{equation}
holds for all $\psi^*\in H^{-1/2}(\Gamma).$ 
The BEM discretization then reads:

Find $\psi_h \in P_h^0\subset H^{-1/2}(\Gamma)$ such that
\begin{equation}\label{eq:SL-weak-discrete}
    \begin{aligned}
    (\text{SL Indirect})\qquad &&\langle \hat{\mathrm{V}}_\omega\psi_h, \psi_h^* \rangle = \langle g_{\mathrm{D}}, \psi_h^*\rangle, 
\end{aligned}
\end{equation}
holds for all $\psi_h^*\in P_h^0 \subset H^{-1/2}(\Gamma).$ \\

To discretize the direct formulations \eqref{eq:BIEs_ext}, resp.~\eqref{eq:BIEs}, we simply replace $g_{\mathrm{D}}$ by $(-\frac{1}{2}\mathrm{I} + \hat{\mathrm{K}}_\omega )g_{\mathrm{D}} \in H^{1/2}(\Gamma)$ in \eqref{eq:SL-weak-discrete}, resp.~by $(\frac{1}{2}\mathrm{I} + \hat{\mathrm{K}}_\omega )g_{\mathrm{D}} \in H^{1/2}(\Gamma)$. For the double-layer equation \eqref{eq:BIEs-indirect_ext} of the exterior problem, we proceed in an analogous way: 

For given $g_{\mathrm{D}}\in H^{1/2}(\Gamma)$, find $\varphi_h \in P_h^1 \subset H^{1/2}(\Gamma)$ such that
\begin{equation}\label{eq:DL-weak-discrete}
    \begin{aligned}
        (\text{DL Indirect})\qquad &&\langle (\tfrac{1}{2}\mathrm{I} + \hat{\mathrm{K}}_\omega )\varphi_h, \varphi^*_h\rangle = \langle g_{\mathrm{D}}, \varphi^*_h\rangle,
    \end{aligned}
\end{equation}
holds for all $\varphi^*_h\in P_h^1 \subset H^{-1/2}(\Gamma)$. 

In the case of the interior problem \eqref{eq:BIEs-indirect}, the operator $\tfrac{1}{2}\mathrm{I} + \hat{\mathrm{K}}_\omega$ is replaced by $-\tfrac{1}{2}\mathrm{I} + \hat{\mathrm{K}}_\omega$ on the left hand side of \eqref{eq:DL-weak-discrete}. 

Finally, the direct formulations involving $\tfrac{1}{2}\mathrm{I} \pm \hat{\mathrm{K}}_\omega$ in \eqref{eq:BIEs_ext}, resp.~\eqref{eq:BIEs}, are discretized in $H^{-1/2}(\Gamma)$: 

For given $g_{\mathrm{s}}\in H^{-1/2}(\Gamma)$, find $\psi_h \in P_h^0 \subset H^{-1/2}(\Gamma)$ such that
\begin{equation}
    \begin{aligned}
        (\text{DL direct})\qquad &&\langle (\tfrac{1}{2}\mathrm{I} \pm \hat{\mathrm{K}}_\omega )\psi_h, \psi^*_h\rangle = \langle g_{\mathrm{s}}, \psi^*_h\rangle,
    \end{aligned}
\end{equation}
holds for all $\psi^*_h\in P_h^1 \subset H^{-1/2}(\Gamma)$.

\remark{For all above formulations, the same convergence rates as for the corresponding boundary element formulation of the Laplace equation are expected. In particular, convergence rates are reduced in the presence of geometric singularities like corners or edges, as discussed further in the numerical experiments in Section \ref{sec:numerics}.
We refer to \cite{gwinnerstephan} for a detailed discussion in the case of the Laplace equation.}

\subsection{Algorithmic considerations}

\subsubsection{Approximation of boundary integral operators and potentials}
The decomposition of the fundamental solution leads to a corresponding decomposition of the potential operators and therefore of the entries of the Galerkin matrix. Given a density $\varphi$ and $\nex\in \mathbb{R}^3\setminus \Gamma$, we obtain the single layer potential as
\begin{align*}
    \hat{\mathrm{S}}_\omega\varphi(\nex) &= \int_{\Gamma}\hat{G}_\omega(\nex, \ney)\varphi(\ney)dS_{\ney} \\
    &= \int_{\Gamma} \hat{G}^{\text{sing}}_{\omega}(\nex, \ney)\varphi(\ney)\dsy +  \int_{\Gamma} \int_0^{\tau_0}\exp(-i\omega\tau)\Gdiff(\nex, \ney, \tau)\varphi(\ney)d\tau\ \dsy \\
    &\qquad+
    \int_{\Gamma} \int_{\tau_0}^{\infty}\exp(-i\omega \tau)\, G(\nex, \ney, \tau)\varphi(\ney)d\tau\ \dsy \\
    &\coloneq \Sw^{\text{sing}}\varphi(\nex) + \Sw^{\text{diff}}\varphi(\nex) + \Sw^{\text{tail}}\varphi(\nex).
\end{align*}
The single layer operator $\hat{\mathrm{V}}_{\omega}$, where $\nex\in \Gamma$, is decomposed in the same way, and for the double layer potential we have for a given density $\varphi$ and $\nex\in \mathbb{R}^3\setminus \Gamma$:
\begin{align*}
    \hat{\mathrm{D}}_\omega\varphi(\nex) &= \int_{\Gamma} \dfrac{\partial\hat{G}_\omega}{\partial{n_{\ney}}}(\nex, \ney)\varphi(\ney)\dsy \\
    &= \int_{\Gamma} \dfrac{\partial\hat{G}_{\text{sing},\omega} }{\partial{n_{\ney}}} (\nex, \ney)\varphi(\ney)\dsy +  \int_{\Gamma}  \int_0^{\tau_0}\exp(-i\omega \tau)\,\dfrac{\partial \Gdiff}{\partial{n_{\ney}}}(\nex, \ney, \tau)\varphi(\ney)d\tau\ \dsy \\
    &\qquad +
    \int_{\Gamma} \int_{\tau_0}^{\infty}\exp(-i\omega \tau)\, \dfrac{\partial G}{\partial{n_{\ney}}} (\nex, \ney, \tau)\varphi(\ney)d\tau\ \dsy\\
    &\coloneqq \Dw^{\text{sing}}\varphi(\nex) + \Dw^{\text{diff}}\varphi(\nex) + \Dw^{\text{tail}}\varphi(\nex).
\end{align*}
For the double layer one should treat the time integral as in Subsections \ref{sec:Ilong} and \ref{sec:Ilong-freq}, adapted to the arising powers of $\tau$. The double layer operator $\Kw$, where $\nex\in \Gamma,$ is decomposed in the same way.

We now define the boundary element matrices related to the sesquilinear forms in \eqref{eq:SL-weak-discrete} and \eqref{eq:DL-weak-discrete}. For a given mesh $\mathcal{T}_h$, with $h>0$
\begin{equation}\label{eq:matrix-defs}
    \begin{aligned}
        \mathbf{V}_{\omega, h} &\coloneqq \mathbf{V}_{\omega,h}^{\text{sing}} + \mathbf{V}_{\omega,h}^{\text{diff}} + \mathbf{V}_{\omega,h}^{\text{tail}}, \\
        \mathbf{K}_{\omega, h} &\coloneqq \mathbf{K}_{\omega,h}^{\text{sing}} + \mathbf{K}_{\omega,h}^{\text{diff}} + \mathbf{K}_{\omega,h}^{\text{tail}}, \\
    \end{aligned}
\end{equation}
with entries 
\begin{equation}\label{eq:matrix-entries}
    \begin{aligned}
    (\mathbf{V}^{\star}_{\omega, h})_{ij} \coloneqq \langle \Vw^{\star}\psi^{(j)}_h, \psi^{(i)}_h\rangle, \\
    (\mathbf{K}^{\star}_{\omega, h})_{ij} \coloneqq \langle \Kw^{\star}\varphi^{(j)}_h, \varphi^{(i)}_h\rangle, \\
    \end{aligned}
\end{equation}
for $\star = \{\text{sing}, \text{diff}, \text{tail}\}$, and basis functions $\{ \psi^{(\ell)}_h\} \subset P^0_h, \ \{\varphi_h^{(\ell)}\}\subset P^1_h$. Corresponding to the properties of the kernels $\hat{G}_{\omega}^{\text{sing}}, \ \Gdiff, $ and $G$, different quadrature rules are used for the assembly of $\mathbf{V}_{\omega, h}^{\text{sing}}, \neV_{\omega, h}^{\text{diff}}$, and $\neV_{\omega, h}^{\text{tail}}$ (or $\mathbf{K}_{\omega, h}^{\text{sing}}, \neK_{\omega, h}^{\text{diff}}$, and $\neK_{\omega, h}^{\text{tail}}$).\\

The kernel $\hat{G}^{\text{sing}}_{\omega}$ is weakly-singular, and standard techniques for weakly-singular kernels and computation of quantities such as 

\begin{equation}\label{eq:Ventry}
(\mathbf{V}^{\text{sing}}_{\omega, h})_{\ell \ell} = \iint\limits_{E_{\ell}\times E_{\ell}} \hat{G}_{\omega}^{\text{sing}}(\nex, \ney)\psi^{(\ell)}(\nex)\psi^{(\ell)}(\ney)\dsy\dsx 
\end{equation}
are used, where $E_{\ell}\coloneqq \text{supp}(\psi_h^{(\ell)})$. Quadrature rules as in Sauter-Schwab \cite[Chapter~5]{sauter2010boundary} allow the accurate evaluation of \eqref{eq:Ventry}.\\

For $\tau> 0$, $G$ in \eqref{eq:fundamental-time1} is a smooth function of $\nex,\ney$. Therefore, the assembly of $\mathbf{V}^{\text{diff}}_{\omega, h}$ or $\mathbf{V}^{\text{tail}}_{\omega, h}$ has no major difficulties. Standard Gauss quadrature can be used for  
\begin{equation}\label{eq:Ventry}
(\mathbf{V}^{\star}_{\omega, h})_{ij} = \iint\limits_{E_{i}\times E_{j}} \hat{G}_{\omega}^{\star}(\nex, \ney)\psi^{(i)}(\nex)\psi^{(j)}(\ney)\dsy\dsx, \quad \text{ for }\star\in \{\text{diff}, \text{tail} \}.
\end{equation}

\subsubsection{Further aspects}

The assembly of the Galerkin matrices $\mathbf{V}_{\omega,h}$ and $\mathbf{K}_{\omega,h}$ requires $\mathcal{O}(N^2)$ storage and arithmetic, where $N$ is the number of degrees of freedom on $\Gamma$. Hierarchical matrix ($\mathcal{H}$-matrix) compression techniques~\cite{bebendorf2008hmatrices, hackbusch2015hierarchical} are directly applicable in our setting, since the kernels      
$\hat{G}_\omega^{\mathrm{diff}}$ and $\hat{G}_\omega^{\mathrm{tail}}$ are smooth away from the diagonal, and hence the corresponding off-diagonal blocks of the Galerkin matrix can be accurately approximated
by low-rank matrices. This reduces storage and matrix-vector complexity to $\mathcal{O}(N \log N)$, making BEM competitive for large-scale simulations.

\section{Numerical experiments}\label{sec:numerics} We study the numerical approximation of the Smoluchowski equation \eqref{eq:smoluchowski_equation_Psi-freq} using its boundary integral formulation both for the  stationary problem ($\omega=0$, in Section \ref{sec:num-stationary}) and for the problem in the frequency domain ($\omega \neq 0$, in Section \ref{sec:num-freq}). 
We look for solutions of \eqref{eq:BIEs} to \eqref{eq:BIEs-indirect} in the appropriate boundary element spaces (see \eqref{eq:bem_spaces}): $$\psi_h,\lambda_h \in P^0_h, \quad \text{ and } \quad \varphi_h, \phi_h \in P^1_h.$$
After solving the relevant boundary integral formulation from \eqref{eq:BIEs} to \eqref{eq:BIEs-indirect}, a postprocessing of solutions is carried out to obtain point evaluations of the solution of the PDE, using \eqref{eq:repformula-exterior}, \eqref{eq:repformula-interior}, resp.~\eqref{eq:rep-formula-indirect}.

Errors of the numerical approximations are computed both for point evaluations and in appropriate Sobolev norms. Given $K$ evaluation points inside the domain, we measure the point evaluation error of the numerical solution $\hat{\Psi}_{\omega,h}$ using the euclidean norm in $\mathbb{R}^{K}$:
\begin{equation}\label{eq:error-1}
 \mathrm{error}(\hat{\Psi}_\omega, \hat{\Psi}_{\omega,h}) \coloneq \dfrac{\left(\displaystyle\sum\limits_{\ell=1}^{K} |\hat{\Psi}_\omega(\mathbf{x_{\ell}}) - \hat{\Psi}_{\omega,h}(\mathbf{x}_{\ell})|^2\right)^{1/2}}{\left(\displaystyle\sum\limits_{\ell=1}^{K} |\hat{\Psi}_\omega(\mathbf{x}_{\ell})|^2\right)^{1/2}}.
\end{equation}
The errors of the densities $\lambda_h$, resp.~$\varphi_h$, in the appropriate Sobolev norms on $\Gamma$ are measured using the equivalent norms given by the discretized elliptic Helmholtz weakly singular operators $\mathbb{V}, \mathbb{V}^{-1}$. Here, $\mathbb{V}$ denotes the weakly singular operator based on the fundamental solution for $(-\Delta + 1)$, given by $\frac{\exp(-|\nex-\ney|)}{4\pi|\nex-\ney|}$.  Then 
\begin{equation}\label{eq:error-2}
\mathrm{error}(\lambda_{\mathrm{ref}},\lambda_h) \coloneq \dfrac{ \langle\lambda_h - \lambda_{\mathrm{ref}}, \mathbb{V}(\lambda_h - \lambda_{\mathrm{ref}})\rangle_{\Gamma}  ^{1/2}}{\langle \lambda_{\mathrm{ref}}, \mathbb{V} \lambda_{\mathrm{ref}}\rangle_{\Gamma}^{1/2}}, \qquad \text{ for }\lambda_h,\lambda_{\mathrm{ref}}\in H^{-1/2}(\Gamma),
\end{equation}
and 
\begin{equation}\label{eq:error-3}
\mathrm{error}(\varphi_{\mathrm{ref}},\varphi_h) \coloneq \dfrac{ \langle\varphi_h - \varphi_{\mathrm{ref}}, \mathbb{V}^{-1}(\varphi_h - \varphi_{\mathrm{ref}})\rangle_{\Gamma}  ^{1/2}}{\langle \varphi_{\mathrm{ref}}, \mathbb{V}^{-1} \varphi_{\mathrm{ref}}\rangle_{\Gamma}^{1/2}}, \qquad \text{ for }\varphi_h,\varphi_{\mathrm{ref}}\in H^{1/2}(\Gamma).
\end{equation}
For discrete solutions in $P_h^0$ and $P_h^1$, Galerkin matrices are used to approximate \eqref{eq:error-2} and \eqref{eq:error-3}.

\begin{figure}[t]
    \centering
    \begin{subfigure}[t]{0.45\textwidth}
    \includegraphics[width=0.62\textwidth]{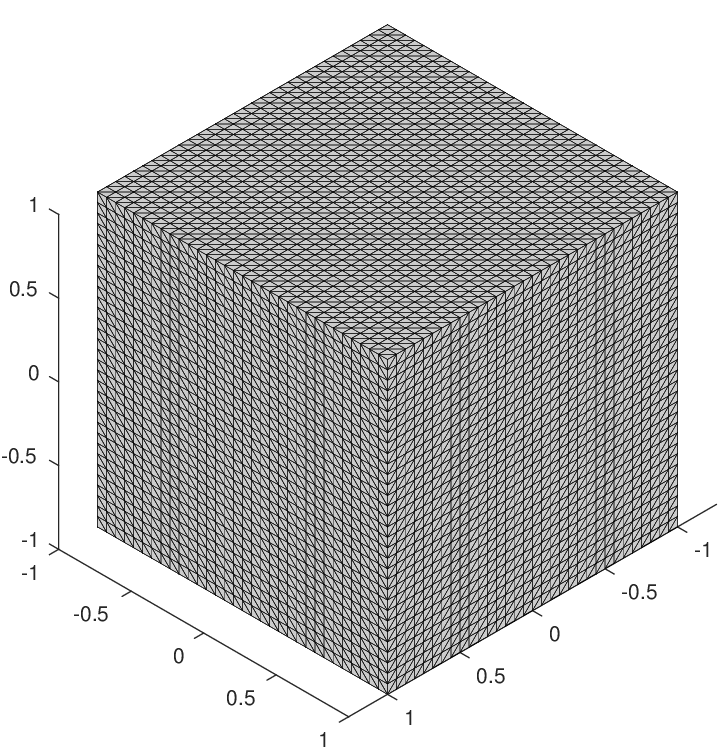}
    \caption{Cube.}
    \label{fig:cube-mesh}
    \end{subfigure}
\begin{subfigure}[t]{0.45\textwidth}
    \centering
    \includegraphics[width=0.75\textwidth]{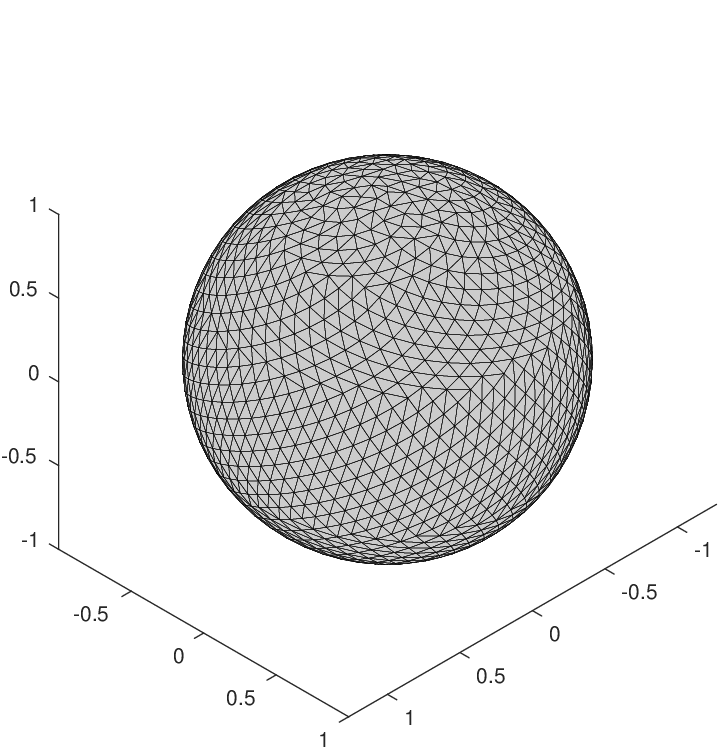}
    \caption{Sphere.}
    \label{fig:sphere-mesh}
\end{subfigure}
\caption{Meshes used in numerical experiments.}
\label{fig:meshes}
\end{figure}

\subsection{Stationary problem ($\omega = 0$)}\label{sec:num-stationary}

\subsubsection{Approximation of fundamental solution: BEM matrices} \label{subsec:FSBEM}

In this experiment we consider the numerical quadrature of the entries in the Galerkin matrices when $\Gamma$ is the boundary of $\Omega_i = [-1,1]^3$ and $\omega=0$. We fix the mesh $\mathcal{T}_h$ depicted in Figure \eqref{fig:cube-mesh} with $12288$ elements.

\subsubsection*{Approximation of $\mathbf{V}^{\text{diff}}_{0, h}$ and $\mathbf{K}^{\text{diff}}_{0, h}$}
First, the matrix entries $(\mathbf{V}^{\text{diff}}_{0, h})_{ij}$ for discretizations in $P_h^0$, resp.~$(\mathbf{K}^{\text{diff}}_{0, h})_{ij}$ for discretizations in $P_h^1$, are studied. The kernel used for the approximation of integral operators $\mathbf{V}^{\text{diff}}_{0, h}$ and $\mathbf{K}^{\text{diff}}_{0, h}$ is continuous in the time variable (see \eqref{eq:continuous}). Therefore, for the accurate approximation we proceed with a composite Gauss quadrature rule. 

\begin{figure}[t]
    \centering
    \begin{subfigure}[t]{0.45\textwidth}
    \includegraphics[width=0.95\textwidth]{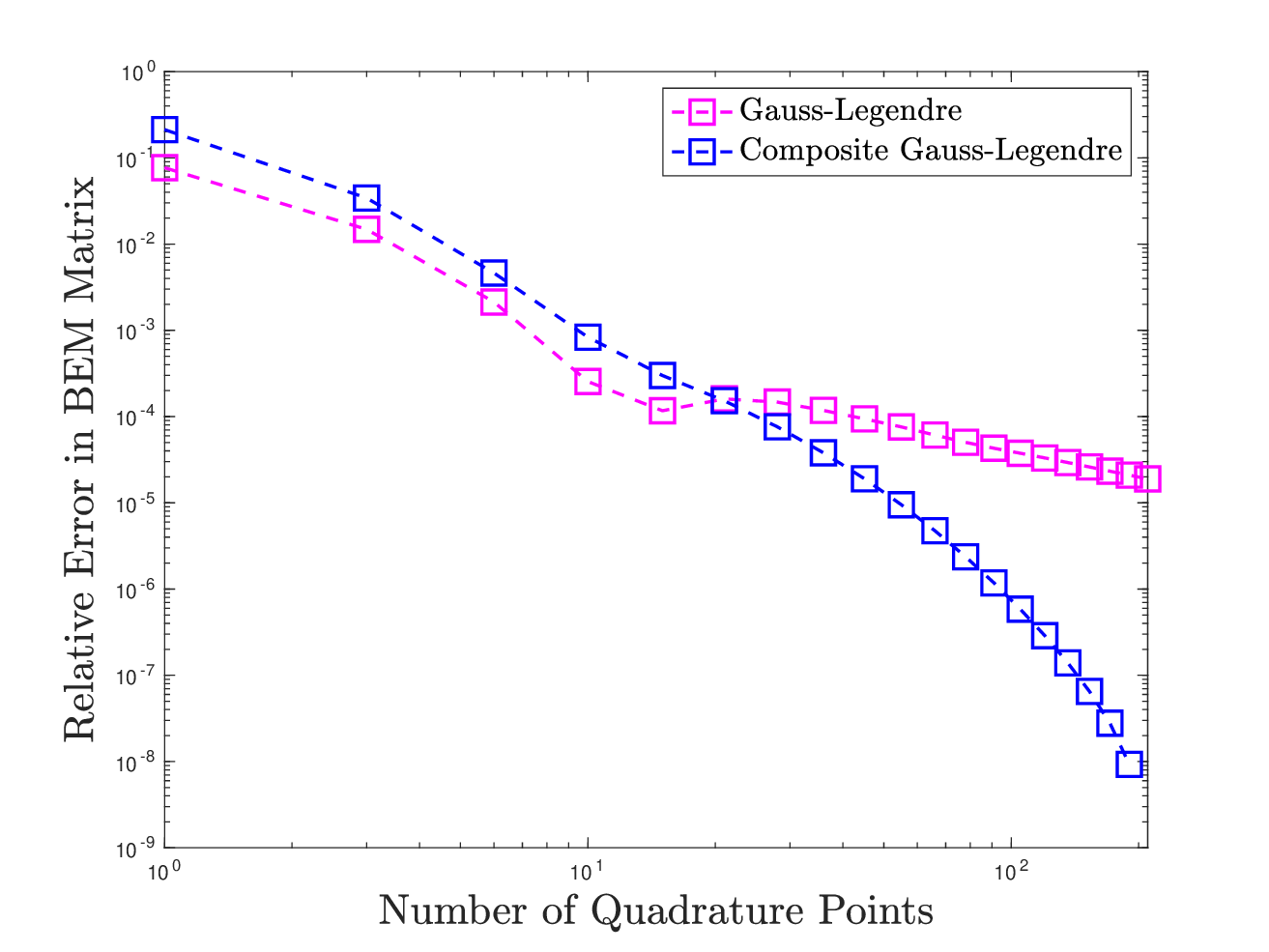}
    \caption{$\mathbf{V}^{\text{diff}}_{0, h}$.}
    \label{fig:Vdiff}
    \end{subfigure}
\begin{subfigure}[t]{0.45\textwidth}
    \centering
    \includegraphics[width=0.95\textwidth]{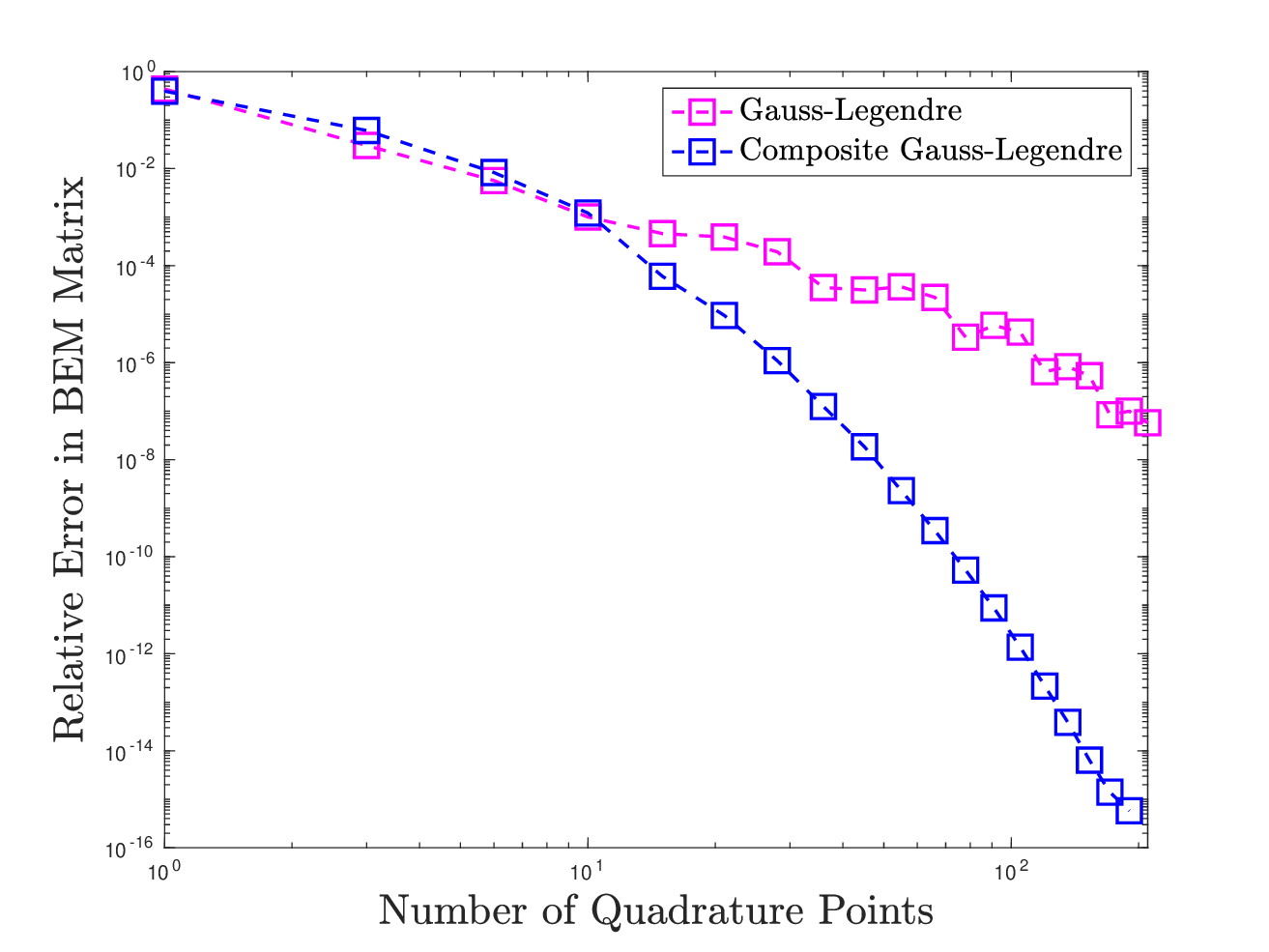}
    \caption{$\mathbf{K}^{\text{diff}}_{0, h}$.}
    \label{fig:Kdiff}
\end{subfigure}
\caption{Error of BEM matrix quadrature for continuous integrand.}
\label{fig:BEM_diff}
\end{figure}

The error of the matrix is measured in the Frobenius (=$2$) norm. Convergence results are shown in Figure \ref{fig:BEM_diff}. Here, we compare with the benchmark obtained from $210$ quadrature points for the composite Gauss-Legendre rule. They show the faster than algebraic  convergence of the composite Gauss-Legendre quadrature as the number of quadrature points increases, reaching a relative error of $10^{-8}$ for around $200$ quadrature points in the case of $\mathbf{V}^{\text{diff}}_{0, h}$, and below $10^{-15}$ in the case of $\mathbf{K}^{\text{diff}}_{0, h}$. The Gauss-Legendre quadrature converges algebraically with rates 1, and 3, respectively.

\subsubsection*{Approximation of $\mathbf{V}^{\text{tail}}_{0, h}$ and $\mathbf{K}^{\text{tail}}_{0, h}$}

\begin{figure}[t]
    \centering
    \begin{subfigure}[t]{0.47\textwidth}
    \includegraphics[width=0.95\textwidth]{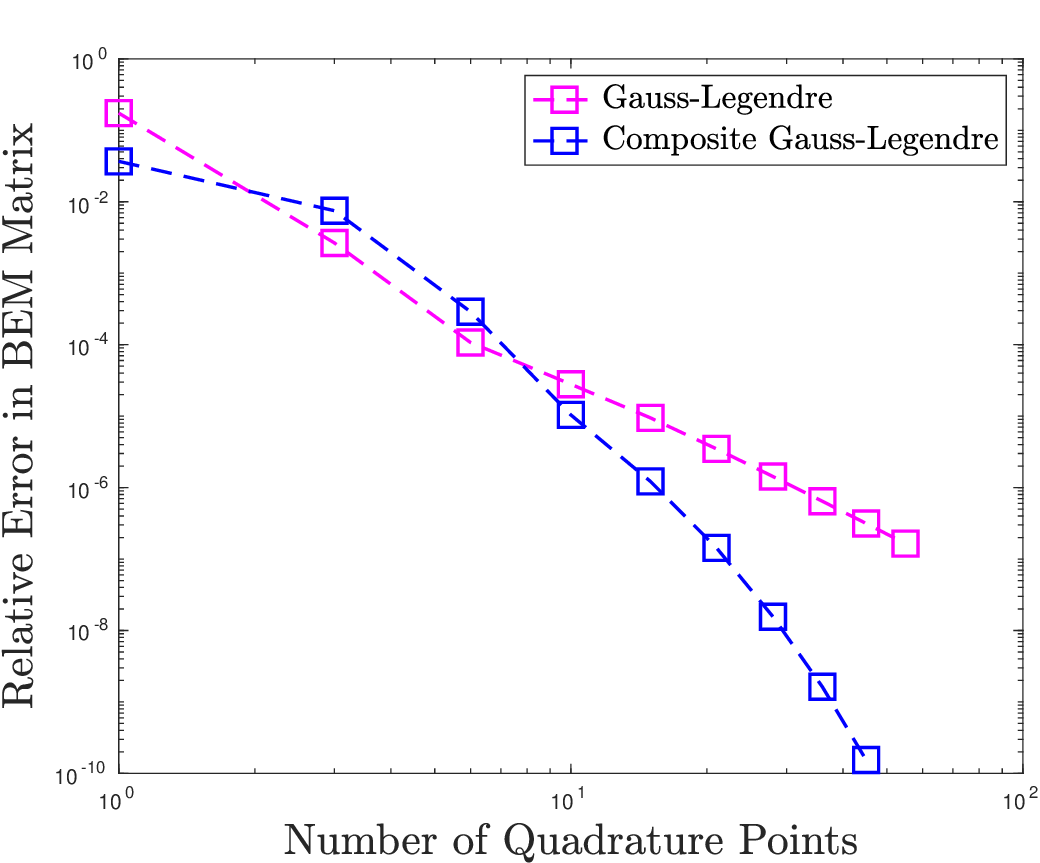}
    \caption{$\mathbf{V}_{0, h}^{\text{tail}}$.}
    \label{fig:Vtail}
    \end{subfigure}
\begin{subfigure}[t]{0.47\textwidth}
    \centering
    \includegraphics[width=0.95\textwidth]{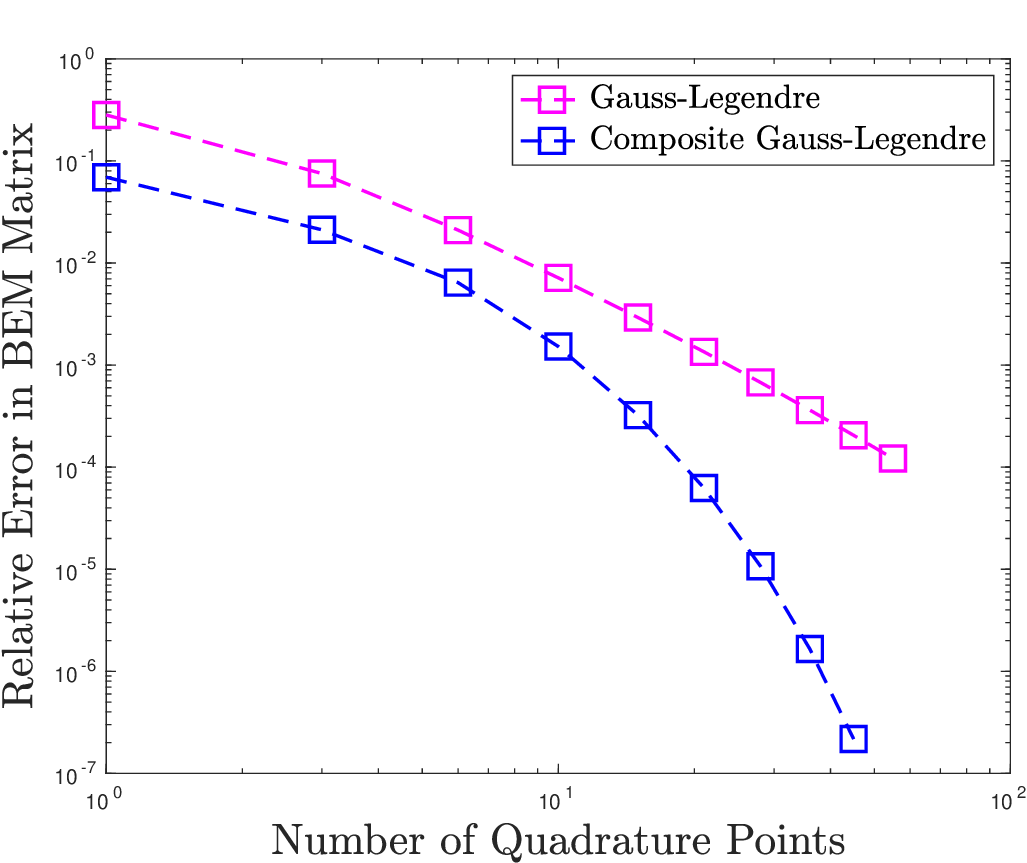}
    \caption{$\mathbf{K}_{0, h}^{\text{tail}}$.}
    \label{fig:Ktail}
\end{subfigure}
\caption{Error of BEM matrix quadrature over unbounded interval. }
\label{fig:BEM_tail}
\end{figure}

We now study the matrix entries $(\mathbf{K}^{\text{tail}}_{0, h})_{ij}$ for discretization in $P_h^0$, resp.~$(\mathbf{K}^{\text{tail}}_{0, h})_{ij}$ for discretization in $P_h^1$. 
The kernel used for the approximation of integral operators $\mathbf{V}^{\text{tail}}_{0, h}$ and $\mathbf{K}^{\text{tail}}_{0, h}$ decays algebraically at infinity in the time variable (see Section \ref{sec:Ilong}). Therefore, the change of variables \eqref{eq:change-variables} reduces the time integral to an integral over a bounded interval with a nonsmooth integrand. Accordingly, a Gauss quadrature leads only to a low order of convergence, and a composite quadrature rule should be used. Convergence results for the matrix error in the Frobenius norm are shown in Figure \ref{fig:BEM_tail}. Here, we compare with the benchmark obtained from $55$ quadrature points for the composite Gauss-Legendre rule. Again, we observe faster than algebraic convergence of the composite Gauss-Legendre quadrature as the number of quadrature points increases, reaching a relative error of $10^{-10}$ for around $45$ quadrature points in the case of $\mathbf{V}^{\text{tail}}_{0, h}$, and $10^{-7}$ in the case of $\mathbf{K}^{\text{tail}}_{0, h}$. The Gauss-Legendre quadrature converges algebraically with rates 3, and 2, respectively. Note that the kernel of $\mathbf{K}^{\text{tail}}_{0, h}$ decays more slowly in time than the kernel of $\mathbf{V}^{\text{tail}}_{0, h}$, explaining the larger errors obtained for $\mathbf{K}^{\text{tail}}_{0, h}$.

\FloatBarrier
\subsubsection{Validation: Interior Problems}\label{sec:validate}

In this experiment we consider the numerical solution of the Smoluchowski equation  inside the domain $\Omega_i = [-1,1]^3$ for $\omega=0$. Boundary conditions on the boundary $\Gamma$ are specified corresponding to the exact solution,
\begin{equation}
\hat{\Psi}_0(x,y,z) = yz + y^2 - z^2 + y^3z - yz^3,
\end{equation}
a harmonic polynomial in the YZ--plane. 

\begin{figure}[htbp]
    \centering
    \begin{subfigure}[b]{0.24\textwidth}
        \centering
        \includegraphics[width=\textwidth]{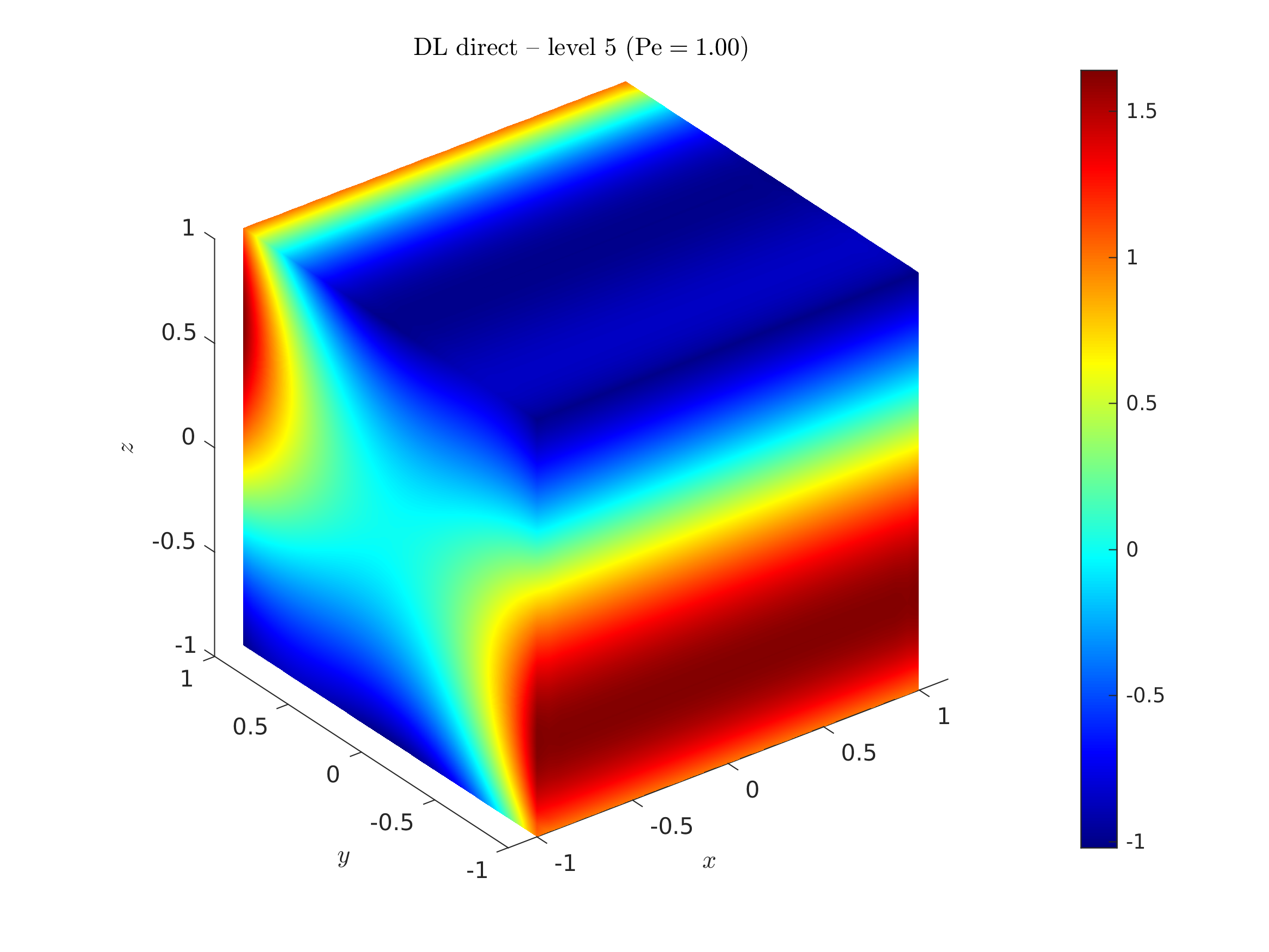}
        \caption{DL direct}
        \label{fig:cube_dl_direct}
    \end{subfigure}
%    \hfill
    \begin{subfigure}[b]{0.24\textwidth}
        \centering
        \includegraphics[width=\textwidth]{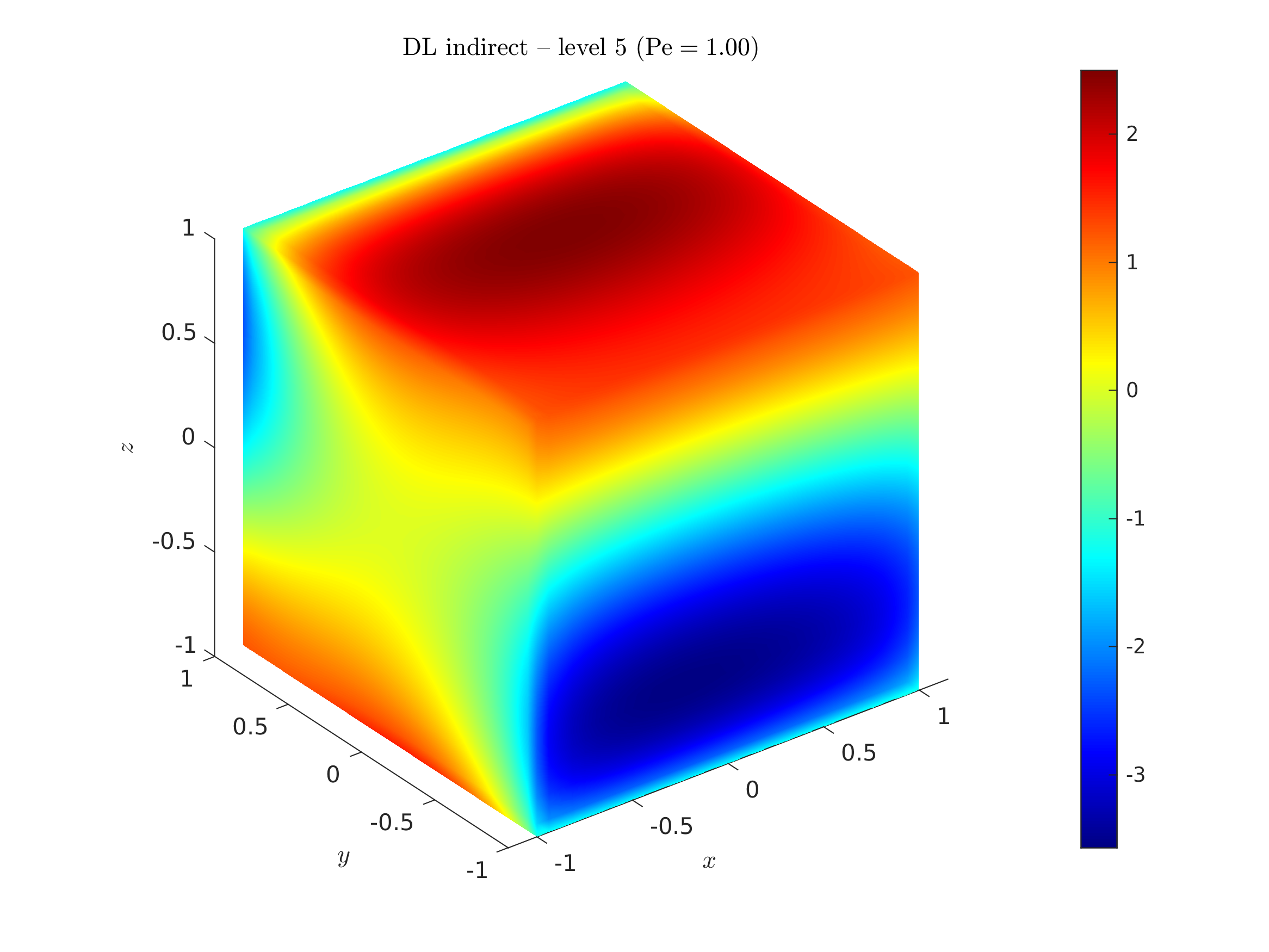}
        \caption{DL indirect}
        \label{fig:cube_dl_indirect}
    \end{subfigure}
%    \vspace{1em}
    \begin{subfigure}[b]{0.24\textwidth}
        \centering
        \includegraphics[width=\textwidth]{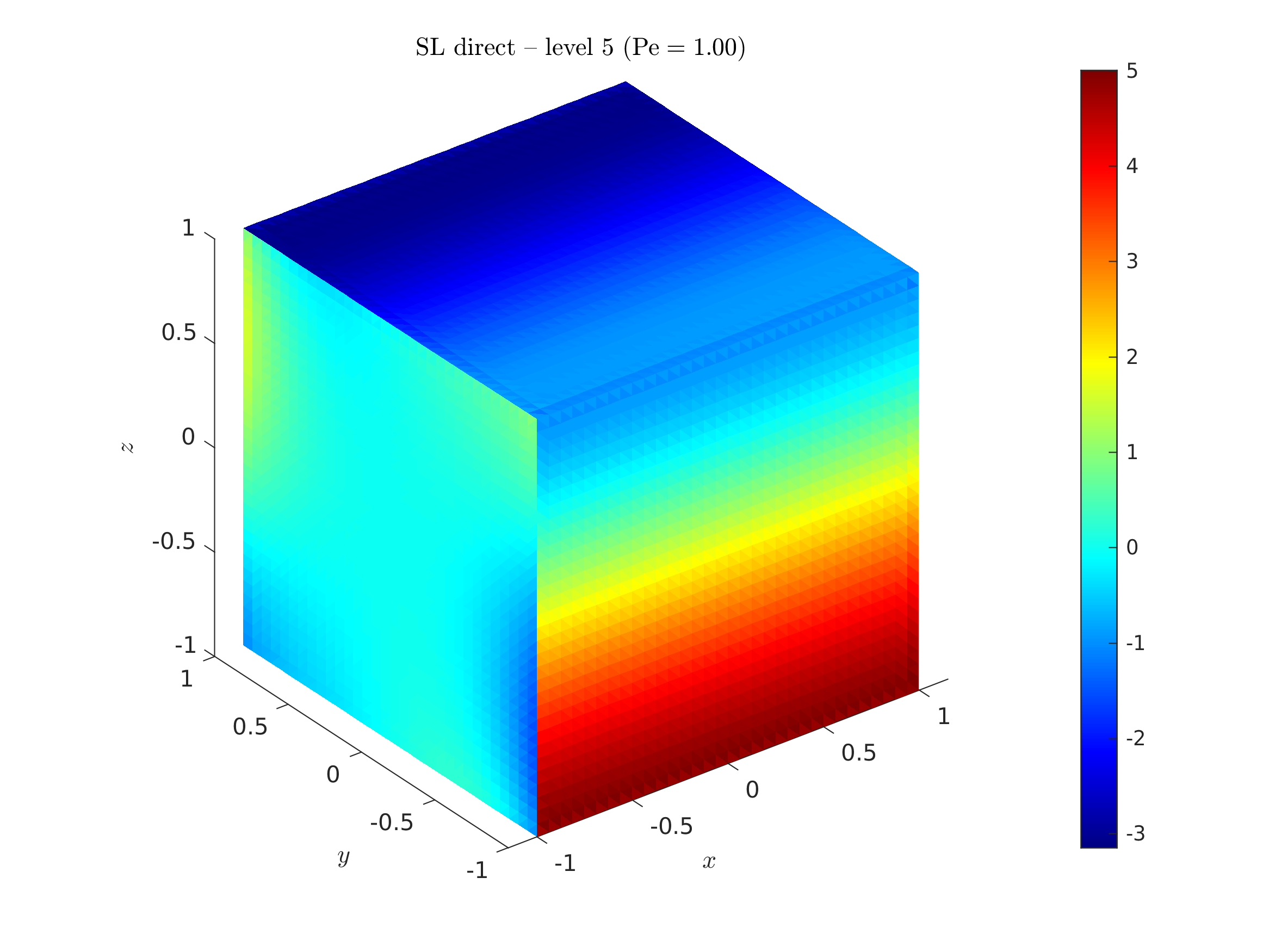}
        \caption{SL direct}
        \label{fig:cube_sl_direct}
    \end{subfigure}
%    \hfill
    \begin{subfigure}[b]{0.24\textwidth}
        \centering
        \includegraphics[width=\textwidth]{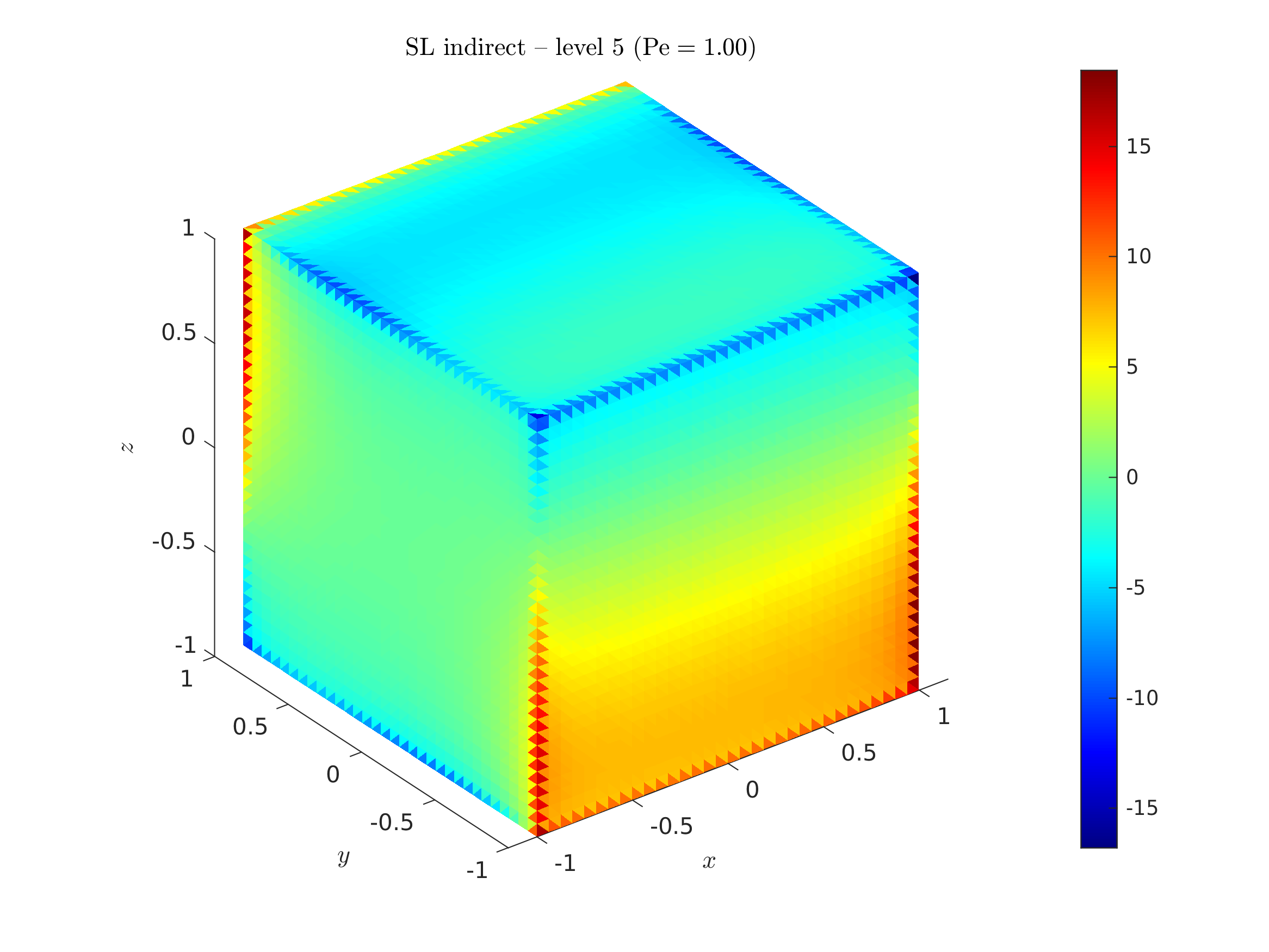}
        \caption{SL indirect}
        \label{fig:cube_sl_indirect}
    \end{subfigure}

    \caption{Snapshots of surface solutions of DL and SL configurations in Section \ref{sec:validate} at the finest mesh, $\Pe = 1$.}
    \label{fig:cube_results}
\end{figure}

Both Dirichlet and Neumann problems are studied, using the
direct and indirect boundary integral equations \eqref{eq:BIEs}, \eqref{eq:BIEs-indirect}. We solve them using the Galerkin discretization from Section \ref{sec:galerkin}, based on the numerically approximated fundamental solution discussed in Section \ref{sec:approx}.
$\Gamma$ is discretized by a sequence of quasi-uniform meshes with up to $12288$ elements. The solutions on the finest mesh are depicted in Figure \ref{fig:cube_results}.

We evaluate the error of point evaluations of the numerical solution to the PDE and the norm error of the density on $\Gamma$, as described at the beginning of this section. To quantify the error of point evaluations, we compare the numerical and exact solutions at $K=10$ points in a circle of radius $0.9$ in the YZ--plane. For the density, we use the solution on the finest mesh as a benchmark.

\begin{figure}[t]
    \centering
    \begin{subfigure}[t]{0.45\textwidth}
    \includegraphics[width=0.95\textwidth]{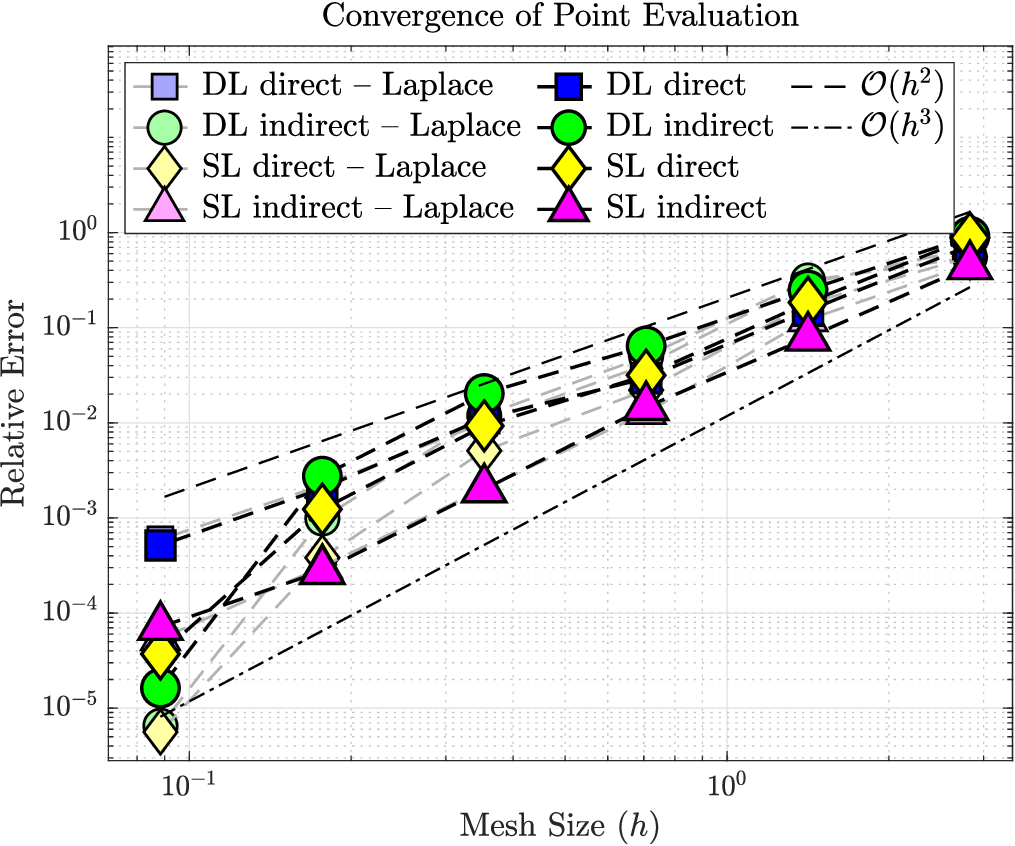}
    \caption{Point evaluations.}
    \label{fig:point-eval-0.25}
    \end{subfigure}
\begin{subfigure}[t]{0.45\textwidth}
    \centering
    \includegraphics[width=0.95\textwidth]{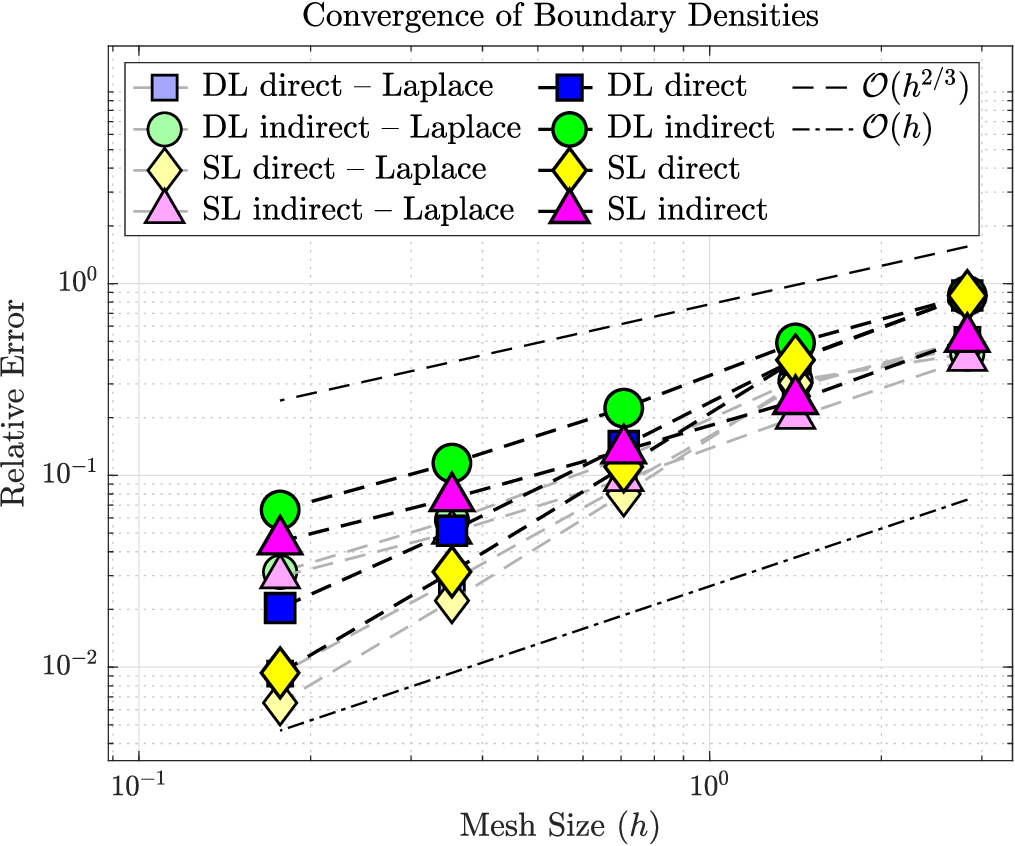}
    \caption{Densities.}
    \label{fig:boundary-0.25}
\end{subfigure}
\caption{Error of solution for interior problem, Section \ref{sec:validate}: $\Pe = \tfrac{1}{4}$. Results for the Laplace equation ($\Pe=0$) are shown in shaded markers for reference.}
\label{fig:conv0.25}
\end{figure}

\begin{figure}[t]
    \centering
    \begin{subfigure}[t]{0.45\textwidth}
    \includegraphics[width=0.95\textwidth]{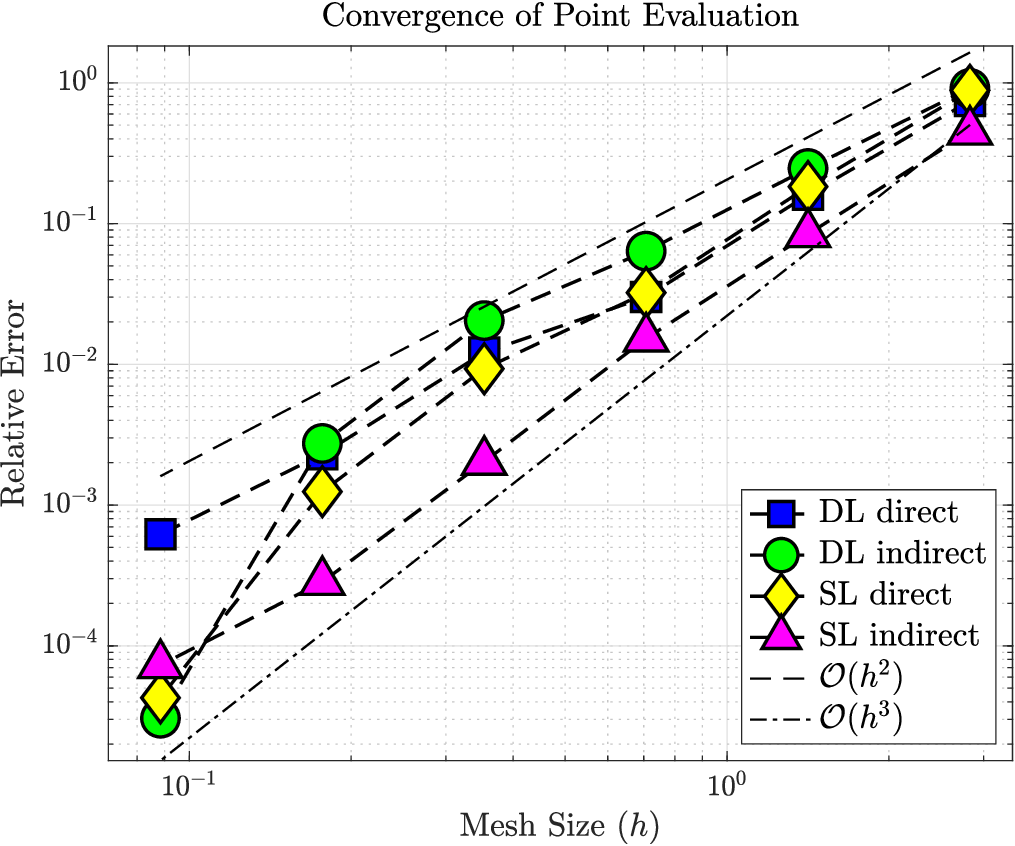}
    \caption{Point evaluations.}
    \label{fig:point-eval-1}
    \end{subfigure}
\begin{subfigure}[t]{0.45\textwidth}
    \centering
    \includegraphics[width=0.95\textwidth]{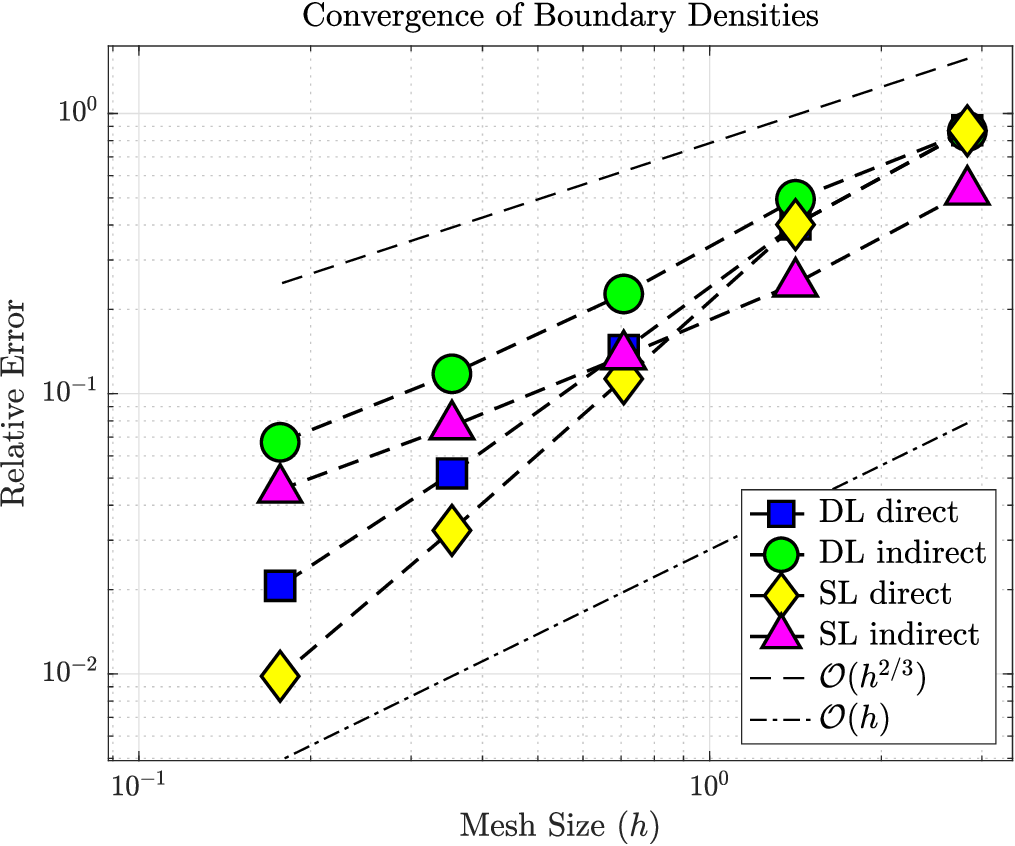}
    \caption{Densities.}
    \label{fig:boundary-1}
\end{subfigure}
\caption{Error of solution for interior problem, Section \ref{sec:validate}: $\Pe = 1$.}
\label{fig:conv1}
\end{figure}

\begin{figure}[t]
    \centering
    \begin{subfigure}[t]{0.45\textwidth}
    \includegraphics[width=0.95\textwidth]{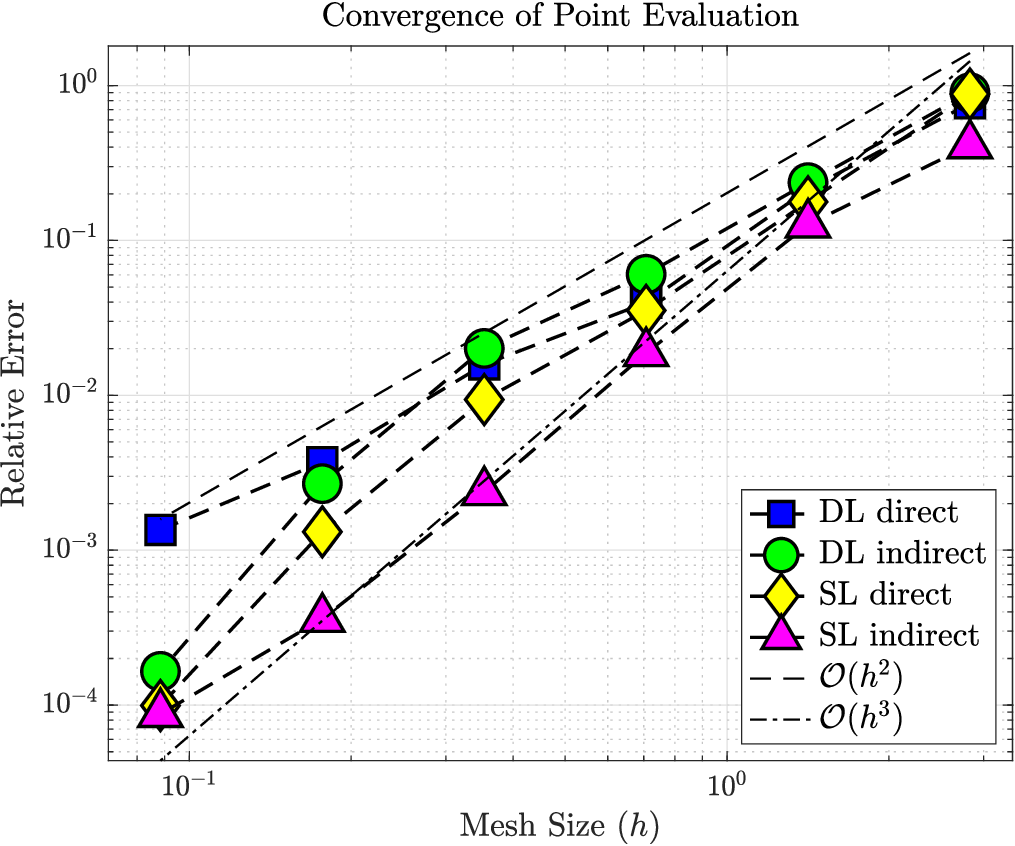}
    \caption{Point evaluations.}
    \label{fig:point-eval-4}
    \end{subfigure}
\begin{subfigure}[t]{0.45\textwidth}
    \centering
    \includegraphics[width=0.95\textwidth]{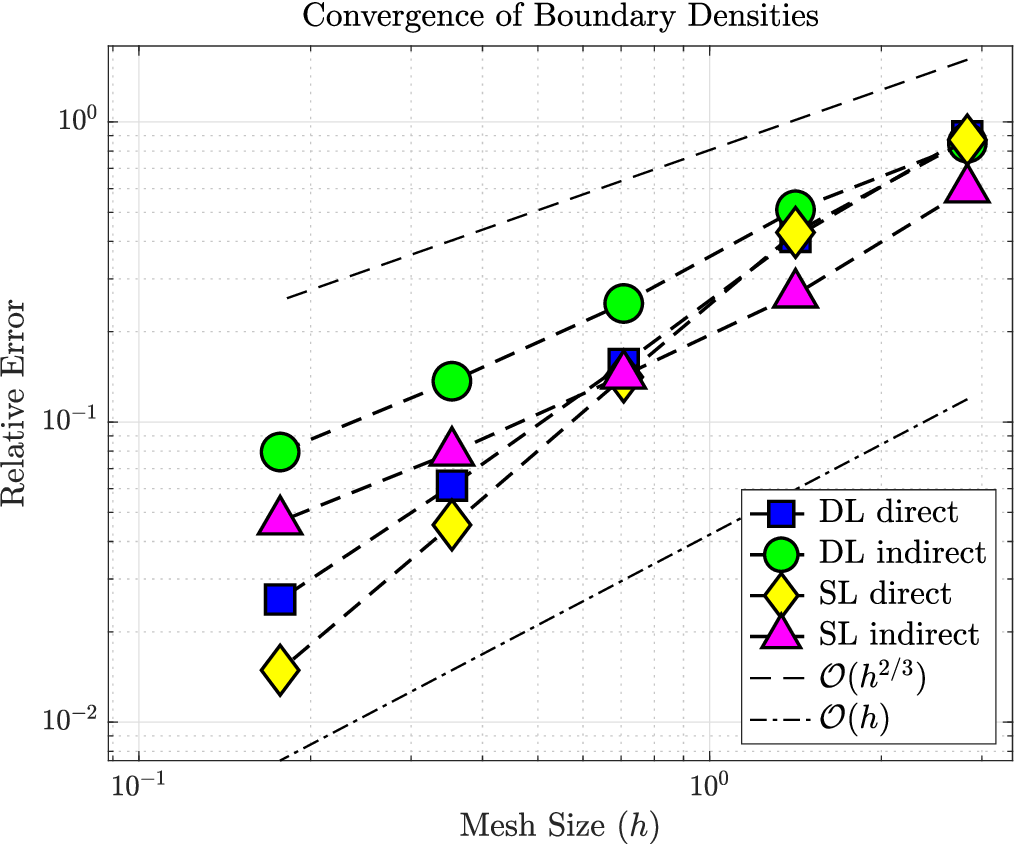}
    \caption{Densities.}
    \label{fig:boundary-4}
\end{subfigure}
\caption{Error of solution for interior problem, Section \ref{sec:validate}: $\Pe = 4$.}
\label{fig:conv4}
\end{figure}

Results in Figures \ref{fig:conv0.25}, \ref{fig:conv1}, and \ref{fig:conv4} are obtained for different P\'eclet numbers, $\Pe \in \{\tfrac{1}{4}, 1, 4\}$. As a reference, we also include convergence results for the Laplacian in Figure \ref{fig:conv0.25}, corresponding to $\Pe = 0$. Similar convergence rates are observed for all problems, including the Laplacian. Due to the nonsmooth boundary of the cube $\Omega$, the indirect formulations have nonsmooth densities as a solution. Correspondingly, the convergence rates are lower than for the direct formulations.

Indeed, edge and corner singularities can be observed in the snapshots of the solutions in Figure \ref{fig:cube_results}, specifically in Figure \ref{fig:cube_sl_indirect} for the indirect single-layer solution and, less pronounced, in Figure \ref{fig:cube_dl_indirect} for the indirect double-layer solution.

\subsection{Frequency domain ($\omega > 0$)}\label{sec:num-freq}

\subsubsection{Approximation of fundamental solution}\label{sec:freq-long}

This subsection studies the accuracy of the numerical quadrature in evaluating the fundamental solution, using the windowing approach presented in Section \ref{sec:Ilong-freq}. We compare the sharp truncation \eqref{eq:integral-long-truncated} and the windowed truncation \eqref{eq:Iwindow} of the integral \eqref{eq:integral-long}, combined with a trapezoidal quadrature rule, resp.~Gauss quadrature. In contrast to \cite{bruno2016windowed, bruno2017windowed}, where the trapezoidal rule is convenient for a periodic integrand, we will observe the advantages of a composite Gauss quadrature rule for \eqref{eq:integral-long}. 

\begin{figure}[t]
    \centering
    \includegraphics[width=0.45\textwidth]{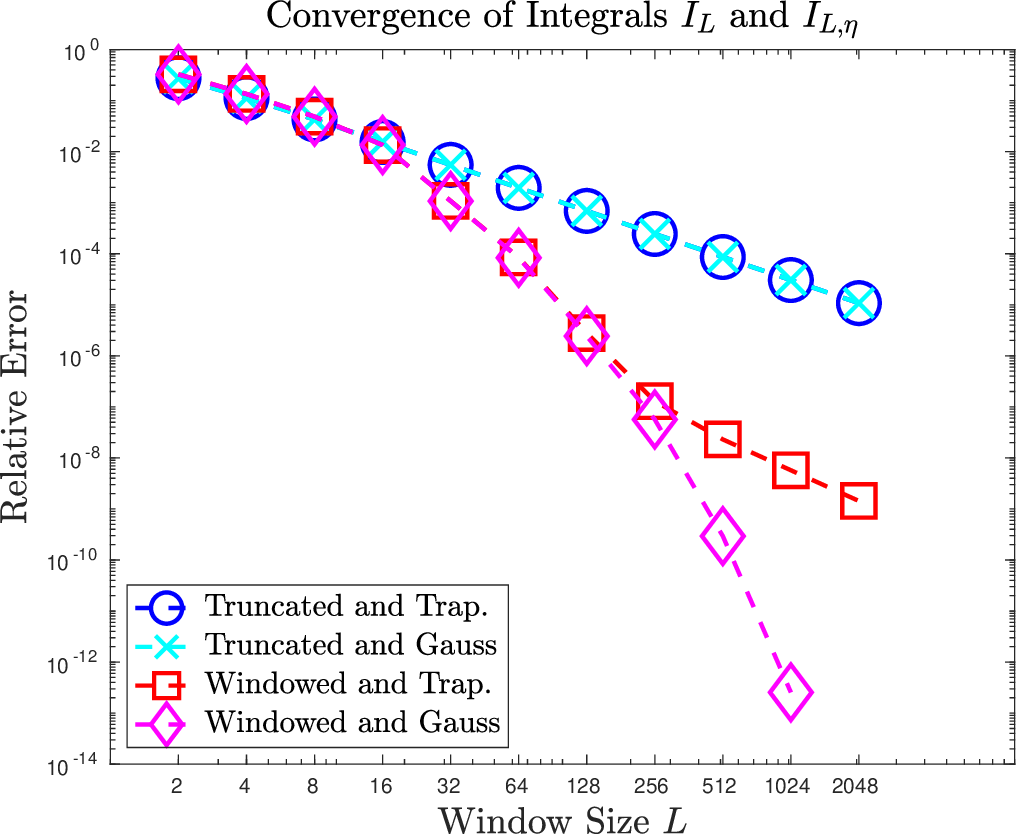}
    \caption{Error of approximations of integral \eqref{eq:integral-long}: sharp, resp.~windowed truncation.}
    \label{fig:window}
\end{figure}

Figure \ref{fig:window} shows the relative error of the numerical quadrature of the integral \eqref{eq:integral-long} with respect to the window size $L$ for $\omega=1$  and $\tau_0=1$. 
For the Gauss quadrature we again use $N_{\textrm{gss}}= 5 \times L$ quadrature points in the interval $[1, L]$. The result for $L = 1024$ is used as a benchmark. We observe faster than algebraic convergence in $L$ for the windowed approximation, until the accuracy of the trapezoidal  ($10^{-7}$), resp.~the Gauss quadrature ($<10^{-12}$) is attained. The truncated approximation, meanwhile, converges algebraically with a rate $\tfrac{3}{2}$. 

\begin{figure}[t]
    \centering
    \begin{subfigure}[t]{0.45\textwidth}
    \includegraphics[width=0.95\textwidth]{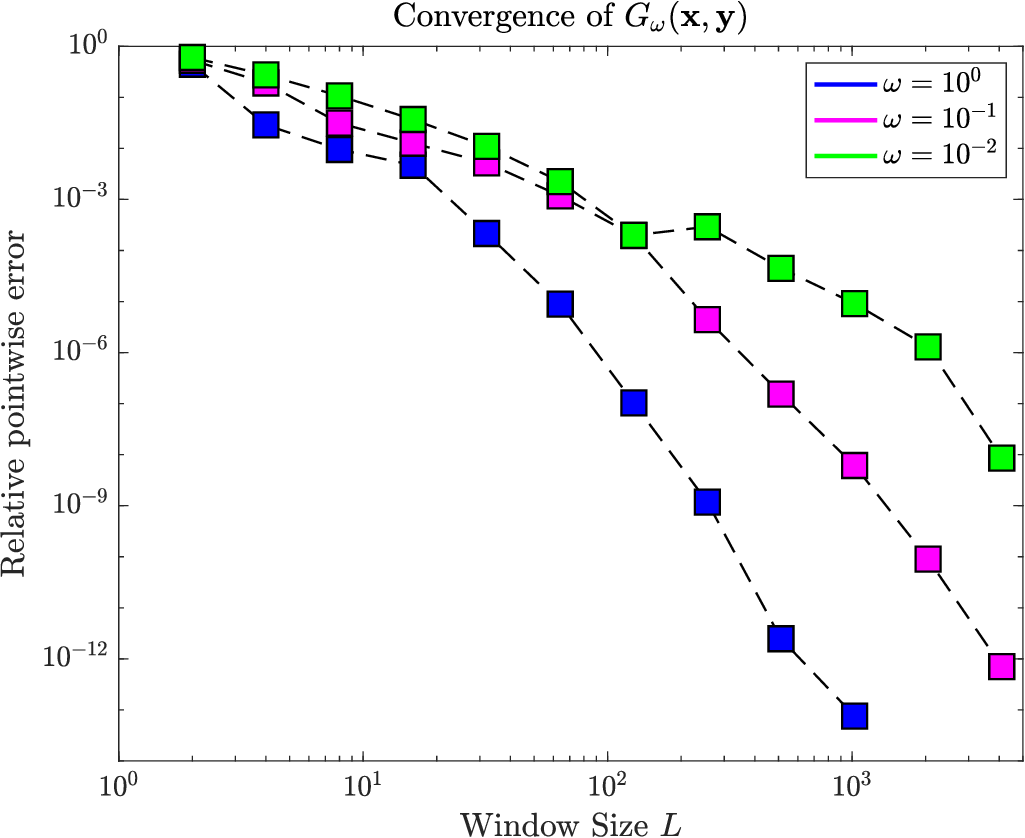}
    \caption{Relative pointwise error of $G$.}
    \label{fig:Gvalues_conv}
    \end{subfigure}
\begin{subfigure}[t]{0.45\textwidth}
    \centering
    \includegraphics[width=0.95\textwidth]{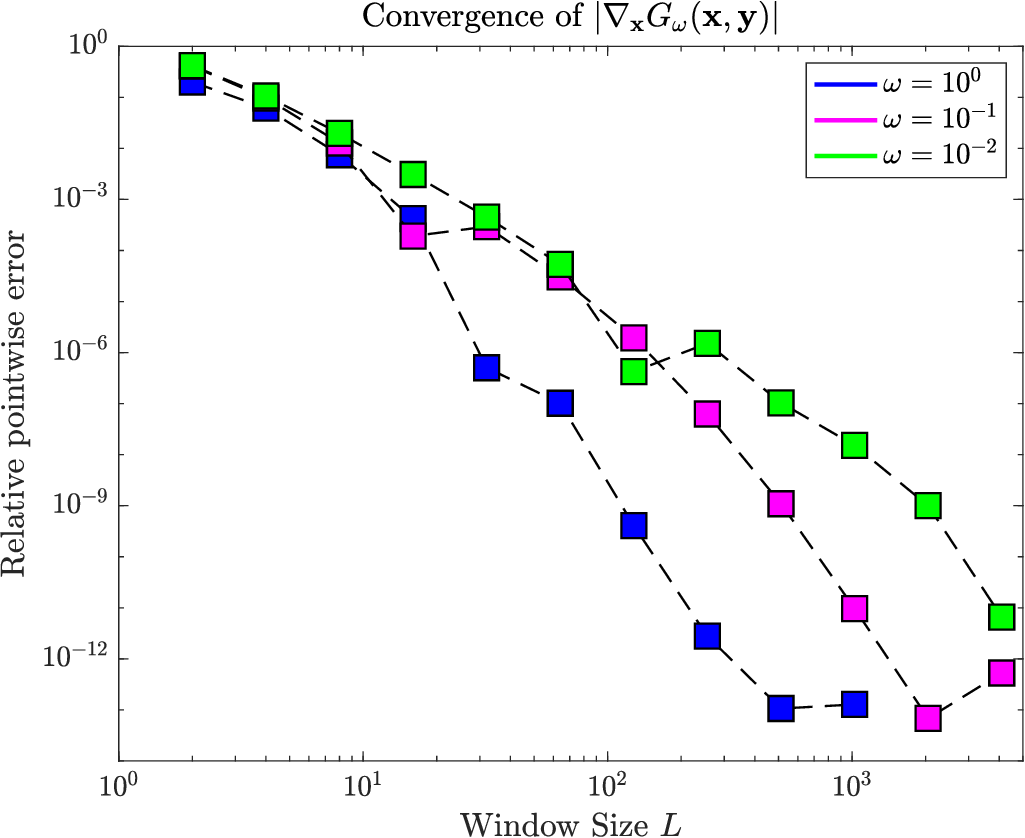}
    \caption{Relative pointwise error of $|\nabla G|$.}
    \label{fig:dGvalues_conv}
\end{subfigure}
\caption{Windowed approximations of $G$ and $\nabla G$ in $\mathbf{x} = (\frac{1}{2}, -\frac{1}{10},\frac{1}{5})$, $\mathbf{y} = (\frac{1}{4}, -\frac{1}{20},\frac{1}{10})$.}
\label{fig:window_conv}
\end{figure}

For the windowed approximation based on the windowing function \eqref{eq:window}, Figure \ref{fig:window_conv} depicts the corresponding relative pointwise error of $G$ and of its gradient $\nabla G$ with respect to the window size $L$. This measures, in particular, the error of the windowed approximation of $\Ilong$. The experiments consider $\Pe = 1$, $\omega\in \{10^{-2}, 10^{-1},10^{0}\}$. Relative errors up to $10^{-13}$ are obtained, where the required $L$ to achieve a given accuracy scales with $\omega$.

% \begin{figure}[t]
%     \centering
%     \includegraphics[width=0.55\textwidth]{figures/laplace-helmholtz.eps}
%     \caption{Convergence of Laplace and Helmholtz approximations of interior problems.}
%     \label{fig:window}
% \end{figure}

\subsubsection{Approximation of fundamental solution: BEM matrices}

We consider the numerical quadrature of the entries in the Galerkin matrices for $\omega = 1$, using the geometry, meshes and postprocessing as in Subsection \ref{subsec:FSBEM}.

\subsubsection*{Approximation of $\mathbf{V}_{\omega, h}^{\text{diff}}$ and $\mathbf{K}_{\omega, h}^{\text{diff}}$}
First, the matrix entries $(\mathbf{V}^{\text{diff}}_{\omega, h})_{ij}$ for discretizations in $P_h^0$, resp.~$(\mathbf{K}^{\text{diff}}_{\omega, h})_{ij}$ for discretizations in $P_h^1$, are studied. As in Subsection \ref{subsec:FSBEM} we use a composite Gauss quadrature rule. 

\begin{figure}[t]
    \centering
    \begin{subfigure}[t]{0.45\textwidth}
    \includegraphics[width=0.95\textwidth]{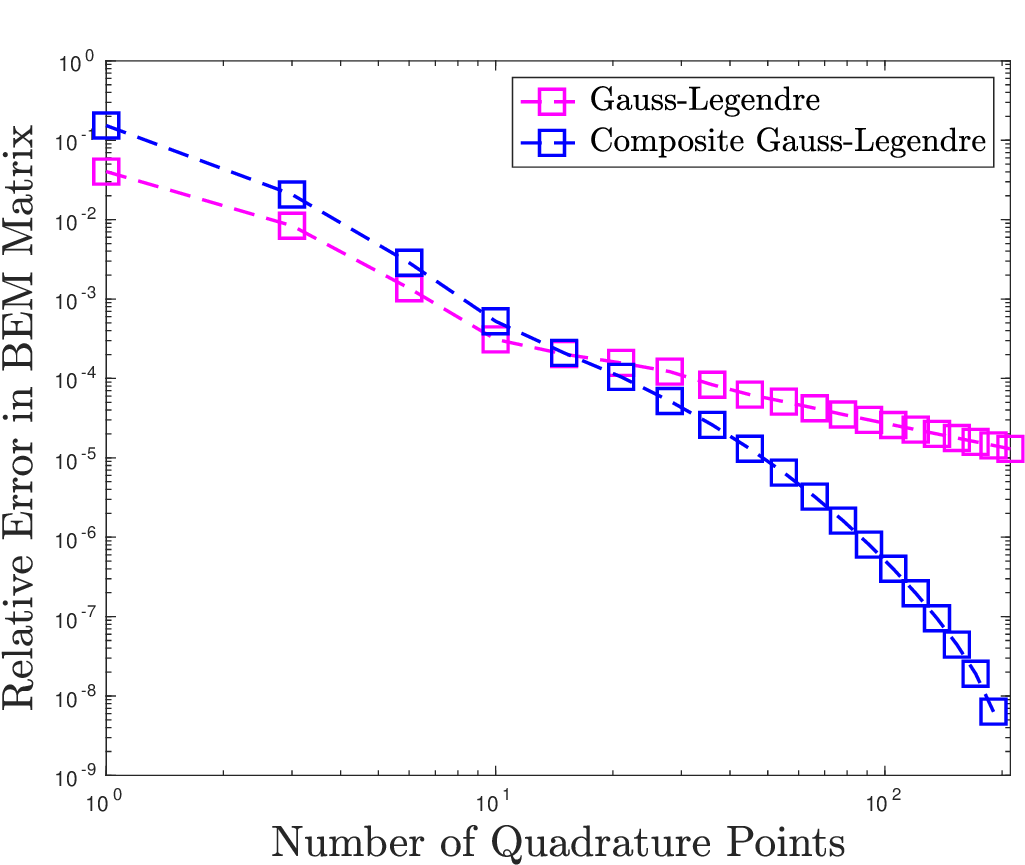}
    \caption{$\mathbf{V}_{\omega, h}^{\text{diff}}$.}
    \label{fig:Vdiff_freq}
    \end{subfigure}
\begin{subfigure}[t]{0.45\textwidth}
    \centering
    \includegraphics[width=0.95\textwidth]{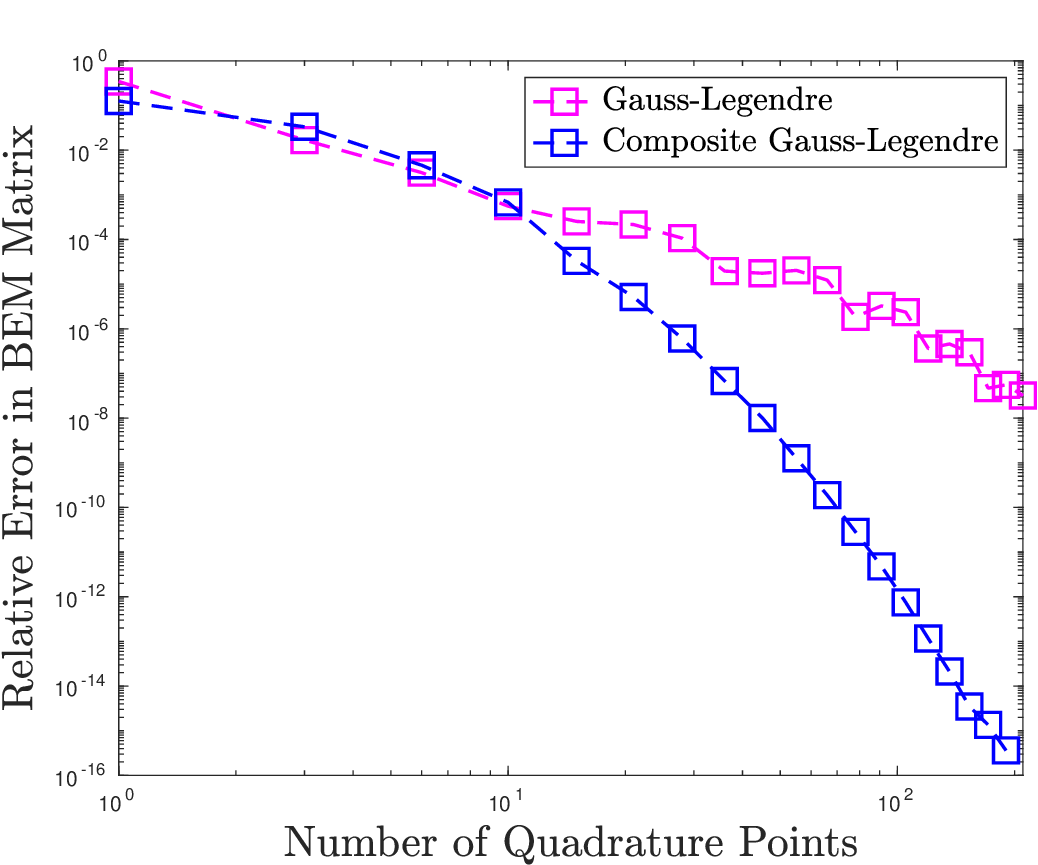}
    \caption{$\mathbf{K}_{\omega, h}^{\text{diff}}$.}
    \label{fig:Kdiff_freq}
\end{subfigure}
\caption{Error of BEM matrix quadrature for continuous integrand ($\omega = 1$).}
\label{fig:BEM_diff_freq}
\end{figure}

The obtained results reflect those for $\omega=0$. Figure \ref{fig:BEM_diff_freq} shows the faster than algebraic  convergence of the composite Gauss-Legendre quadrature as the number of quadrature points increases, reaching a relative error of $10^{-8}$ for around $200$ quadrature points in the case of $\mathbf{V}^{\text{diff}}_{\omega, h}$, and below $10^{-15}$ in the case of $\mathbf{K}^{\text{diff}}_{0, h}$. The Gauss-Legendre quadrature converges algebraically with rates 1 and 3, respectively. 

\subsubsection*{Approximation of $\mathbf{V}_{\omega, h}^{\text{tail}}$ and $\mathbf{K}_{\omega, h}^{\text{tail}}$}
\begin{figure}[t]
    \centering
    \begin{subfigure}[t]{0.47\textwidth}
    \includegraphics[width=0.95\textwidth]{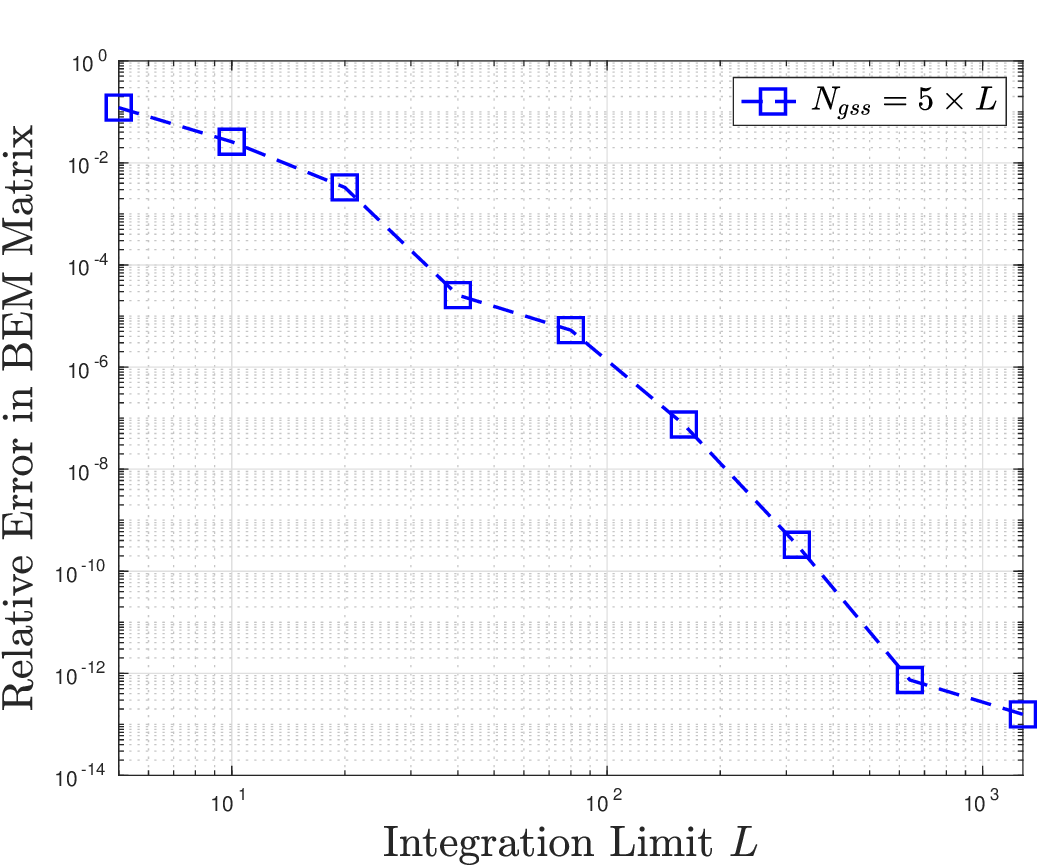}
    \caption{$\mathbf{V}_{\omega, h}^{\text{tail}}$.}
    \label{fig:Vtail_freq}
    \end{subfigure}
\begin{subfigure}[t]{0.47\textwidth}
    \centering
    \includegraphics[width=0.95\textwidth]{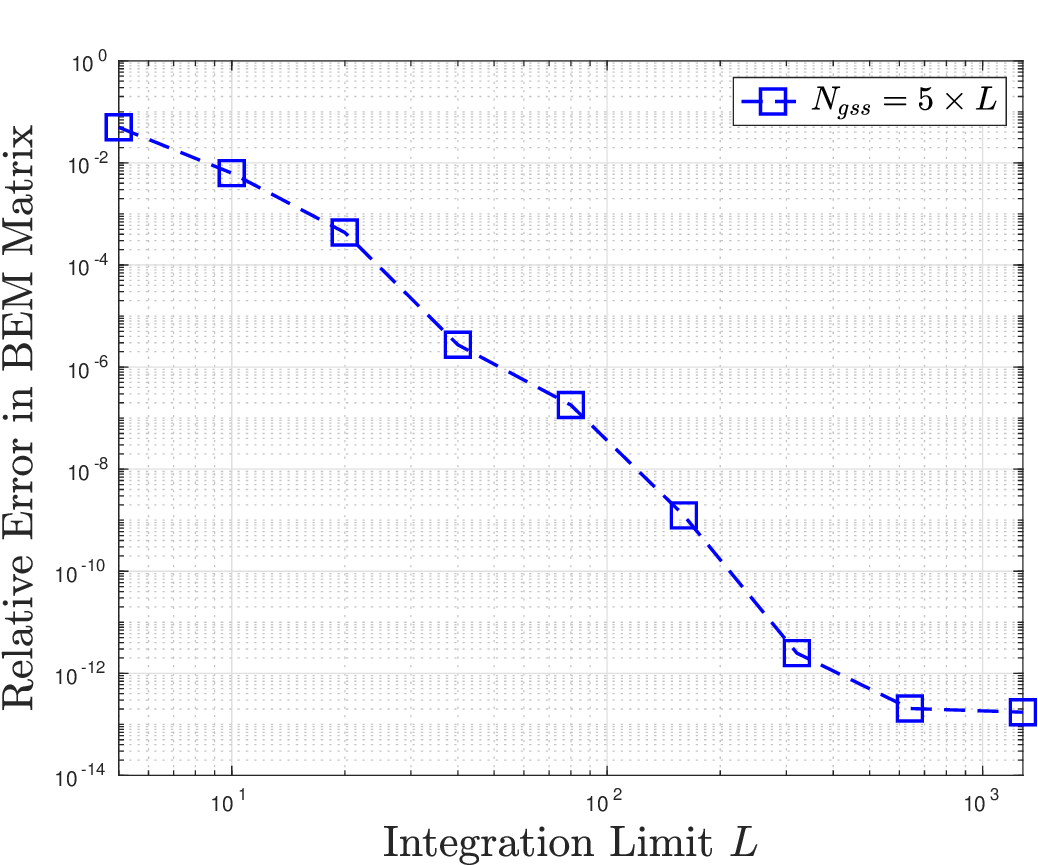}
    \caption{$\mathbf{K}_{\omega, h}^{\text{tail}}$.}
    \label{fig:Ktail_freq}
\end{subfigure}
\caption{Error of BEM matrix quadrature over unbounded interval ($\omega = 1$).}
\label{fig:BEM_tail_freq}
\end{figure}

We now study the matrix entries $(\mathbf{K}^{\text{tail}}_{\omega, h})_{ij}$ for discretizations in $P_h^0$, resp.~$(\mathbf{K}^{\text{tail}}_{\omega, h})_{ij}$ for discretizations in $P_h^1$. 
The windowed kernel used for the approximation of integral operators $\mathbf{V}_{\omega, h}^{\text{tail}}$ and $\mathbf{K}_{\omega, h}^{\text{tail}}$ decays to zero at $\tau = L$ (see Section \ref{sec:Ilong} and \eqref{eq:window}). 
We study the approximation of the Galerkin matrices as the window size $L$ increases, using the approach of Section \ref{sec:freq-long}. Figure \ref{fig:BEM_tail_freq} shows the rapid convergence of the matrix, with errors $<10^{-12}$ for both $\mathbf{V}_{\omega, h}^{\text{tail}}$ and $\mathbf{K}_{\omega, h}^{\text{tail}}$ when $L=512$. Here, results are compared with a benchmark for $L=2048$ and number of Gauss points $N_{gss} = 5 \times L$.

\subsubsection{Validation: Interior Problems}\label{sec:validate-freq}

In this experiment we consider the numerical solution of the Smoluchowski equation  inside the domain $\Omega_i = [-1,1]^3$ for fixed $\Pe=1$ and different $\omega$. Boundary conditions on the boundary $\Gamma$ are specified corresponding to the exact solution given by a plane wave in the YZ--plane in the direction $\mathbf{d} = (0, 1, 0)$,
\begin{equation}
\hat{\Psi}_\omega(\mathbf{x}) = \exp(-\kappa \mathbf{x}\cdot \mathbf{d}).
\end{equation}
Here, $\kappa \coloneq \sqrt{i\omega}=\sqrt{\omega} \,\mathrm{e}^{i\frac{\pi}{4}}$.

The discretization and postprocessing of the direct and indirect boundary integral equations \eqref{eq:BIEs}, \eqref{eq:BIEs-indirect} follows Subsection \ref{sec:validate}, using a sequence of uniform meshes with up to $12288$ elements. 

\begin{figure}[t]
    \centering
    \begin{subfigure}[t]{0.45\textwidth}
    \includegraphics[width=0.95\textwidth]{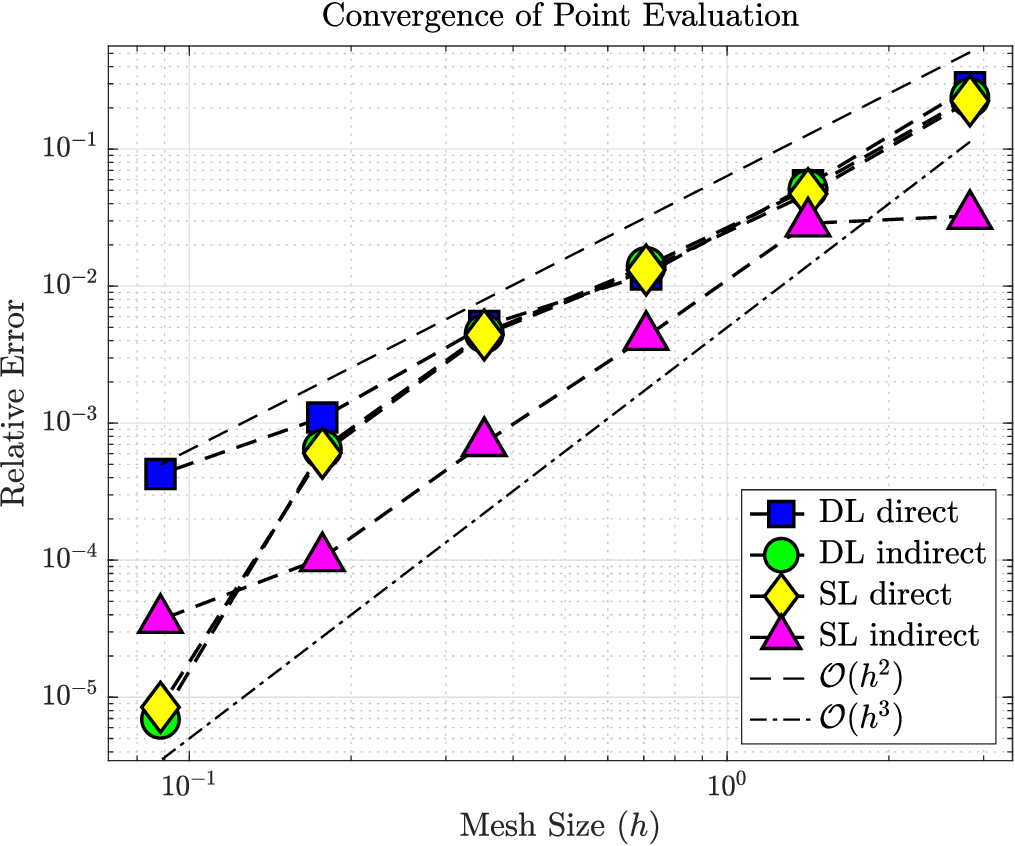}
    \caption{Point evaluations.}
    \label{fig:fpoint-eval-0.25}
    \end{subfigure}
\begin{subfigure}[t]{0.45\textwidth}
    \centering
    \includegraphics[width=0.95\textwidth]{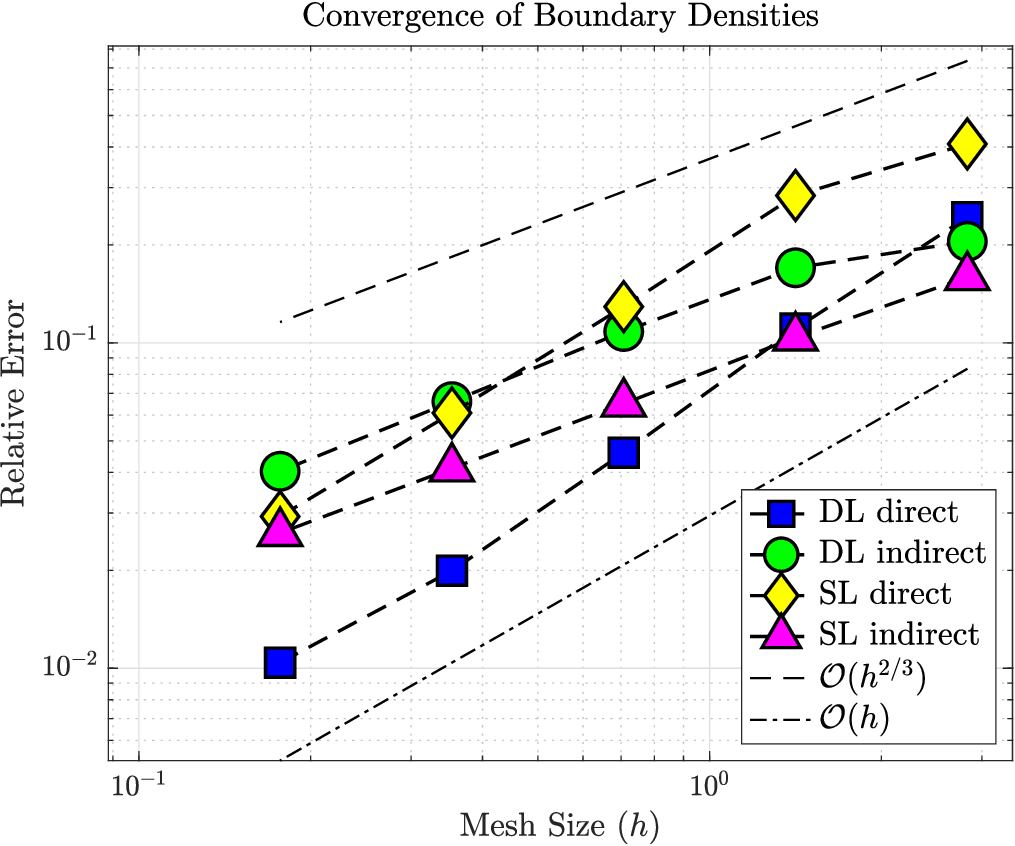}
    \caption{Densities.}
    \label{fig:fboundary-0.25}
\end{subfigure}
\caption{Error of solution for interior problem, Section \ref{sec:validate-freq}: $\omega = \tfrac{1}{4}$.}
\label{fig:fconv0.25}
\end{figure}

\begin{figure}[ht!]
    \centering
    \begin{subfigure}[t]{0.45\textwidth}
    \includegraphics[width=0.95\textwidth]{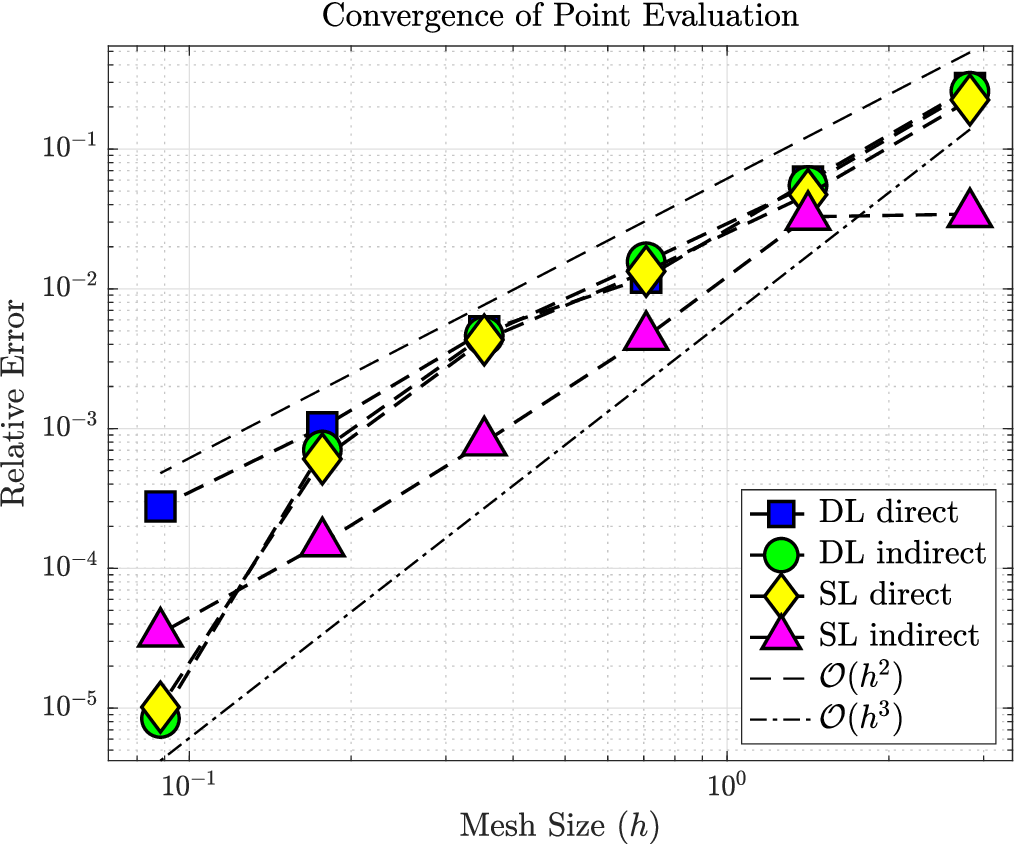}
    \caption{Point evaluations.}
    \label{fig:fpoint-eval-1}
    \end{subfigure}
\begin{subfigure}[t]{0.45\textwidth}
    \centering
    \includegraphics[width=0.95\textwidth]{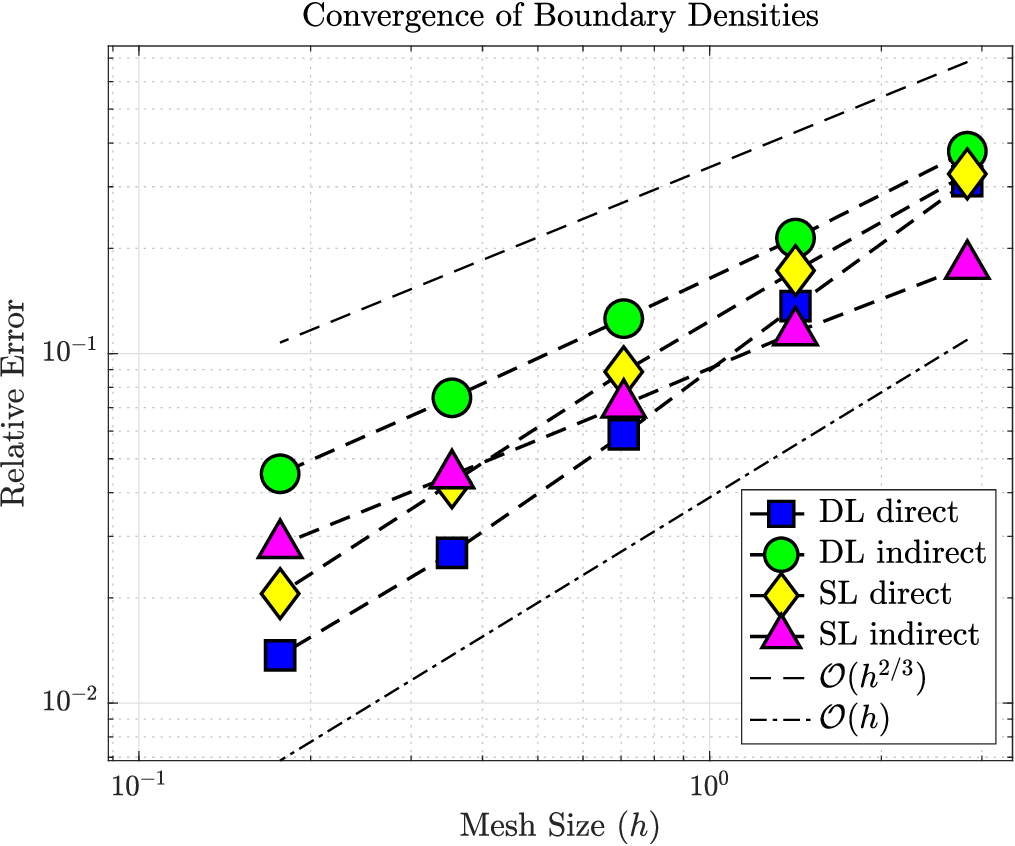}
    \caption{Densities.}
    \label{fig:fboundary-1}
\end{subfigure}
\caption{Error of solution for interior problem, Section \ref{sec:validate-freq}: $\omega = 1$.}
\label{fig:fconv1}
\end{figure}

\begin{figure}[ht!]
    \centering
    \begin{subfigure}[t]{0.45\textwidth}
    \includegraphics[width=0.95\textwidth]{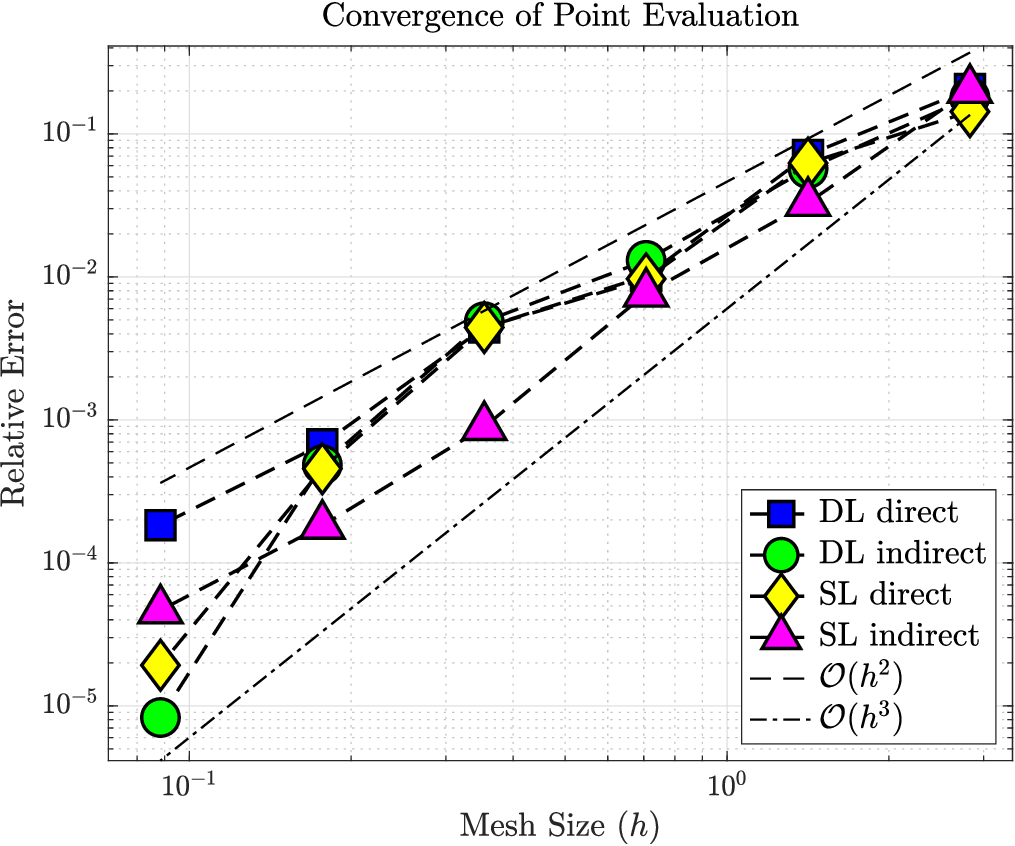}
    \caption{Point evaluations.}
    \label{fig:fpoint-eval-4}
    \end{subfigure}
\begin{subfigure}[t]{0.45\textwidth}
    \centering
    \includegraphics[width=0.95\textwidth]{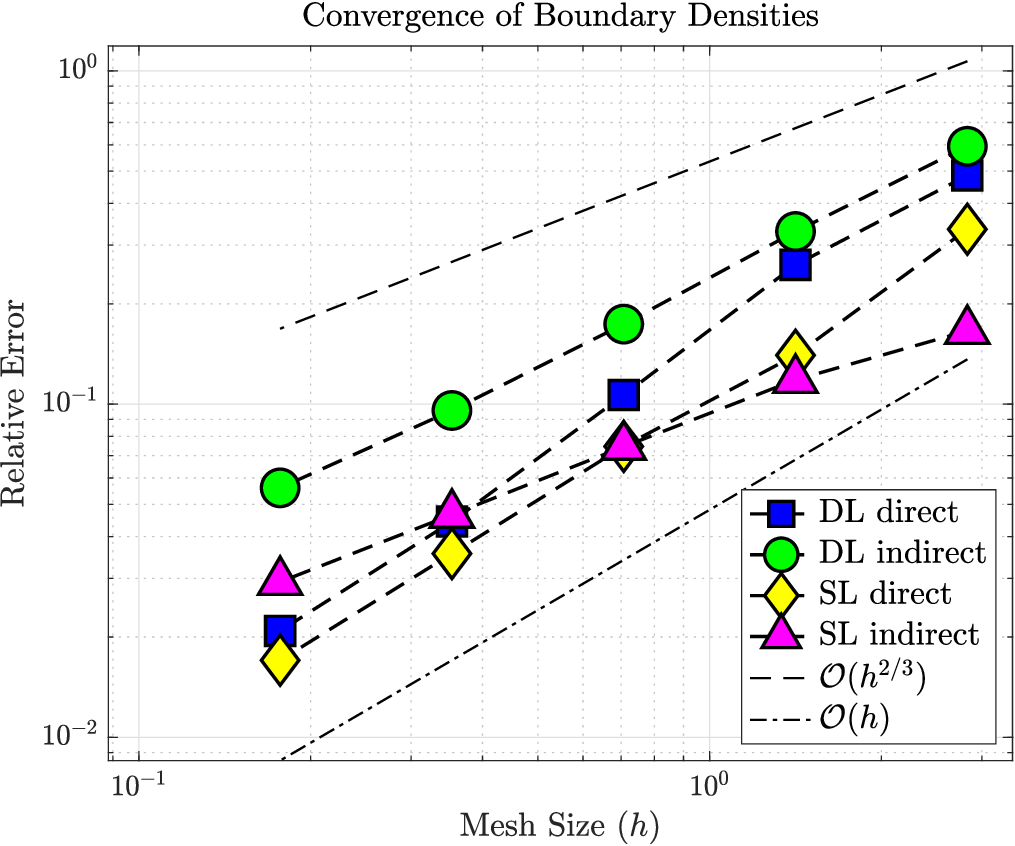}
    \caption{Densities.}
    \label{fig:fboundary-4}
\end{subfigure}
\caption{Error of solution for interior problem, Section \ref{sec:validate-freq}: $\omega = 4$.}
\label{fig:fconv4}
\end{figure}

Figures \ref{fig:fconv0.25}, \ref{fig:fconv1}, and \ref{fig:fconv4} evaluate the error of point evaluations of the numerical solution to the PDE and the norm error of the density on $\Gamma$ for different frequencies $\omega\in \{\tfrac{1}{4}, 1, 4\}$.
Similar convergence rates are observed for all problems, in agreement with those obtained for $\omega=0$ in Subsection \ref{sec:validate}. As described there, the indirect formulations show lower convergence rates than the direct formulations, due to the nonsmooth boundary of the cube $\Omega_i$.

\subsection{Exterior problems for colloidal suspensions}
\begin{figure}[t]
 \begin{subfigure}[t]{0.48\textwidth}

    \centering
    \includegraphics[width=0.75\linewidth]{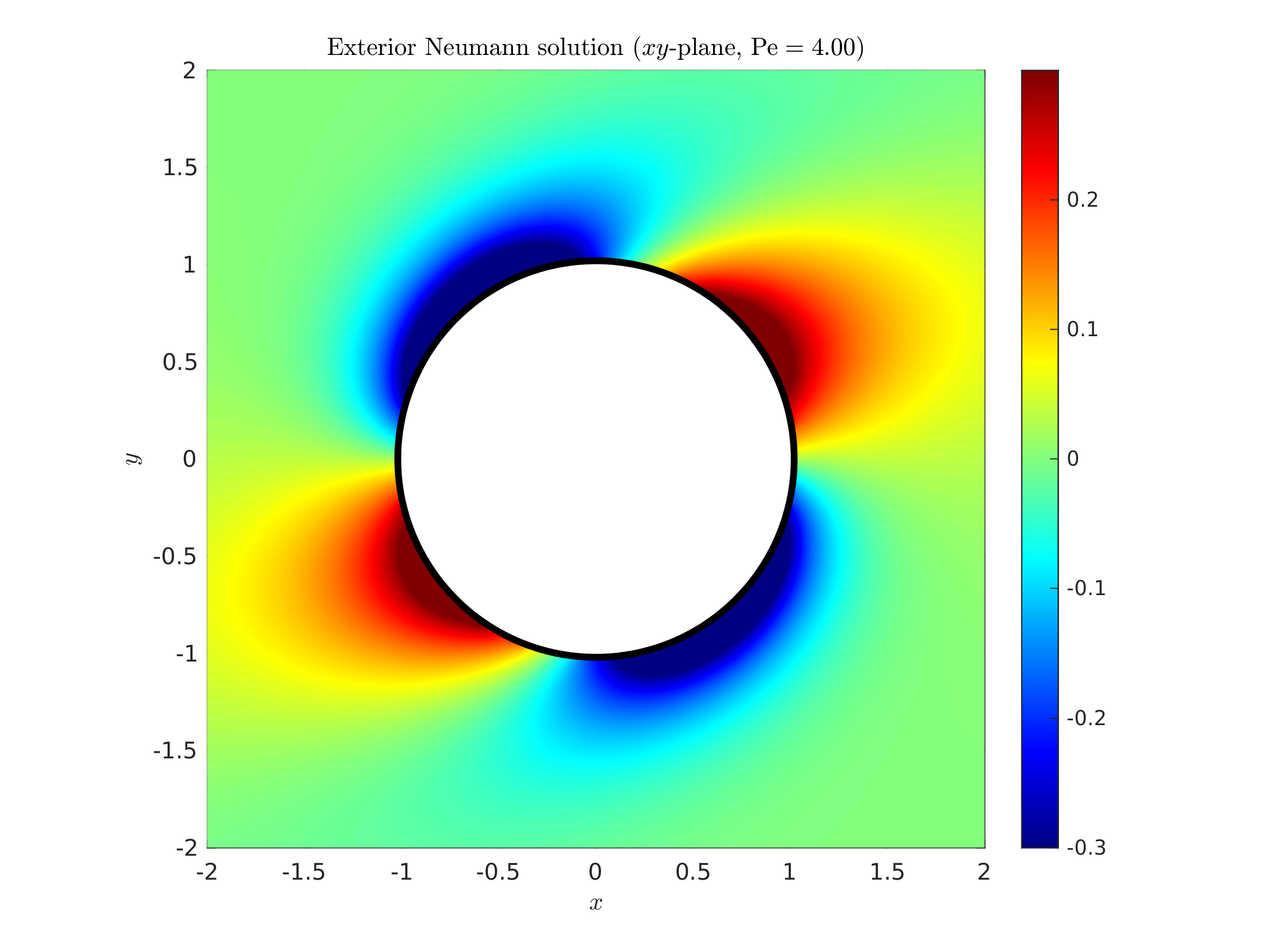}
    \caption{Exterior of one ball.}
    \label{fig:exterior-neumann-sphere}
    \end{subfigure}
 \begin{subfigure}[t]{0.5\textwidth}
     \centering
    \includegraphics[width=0.73\linewidth]{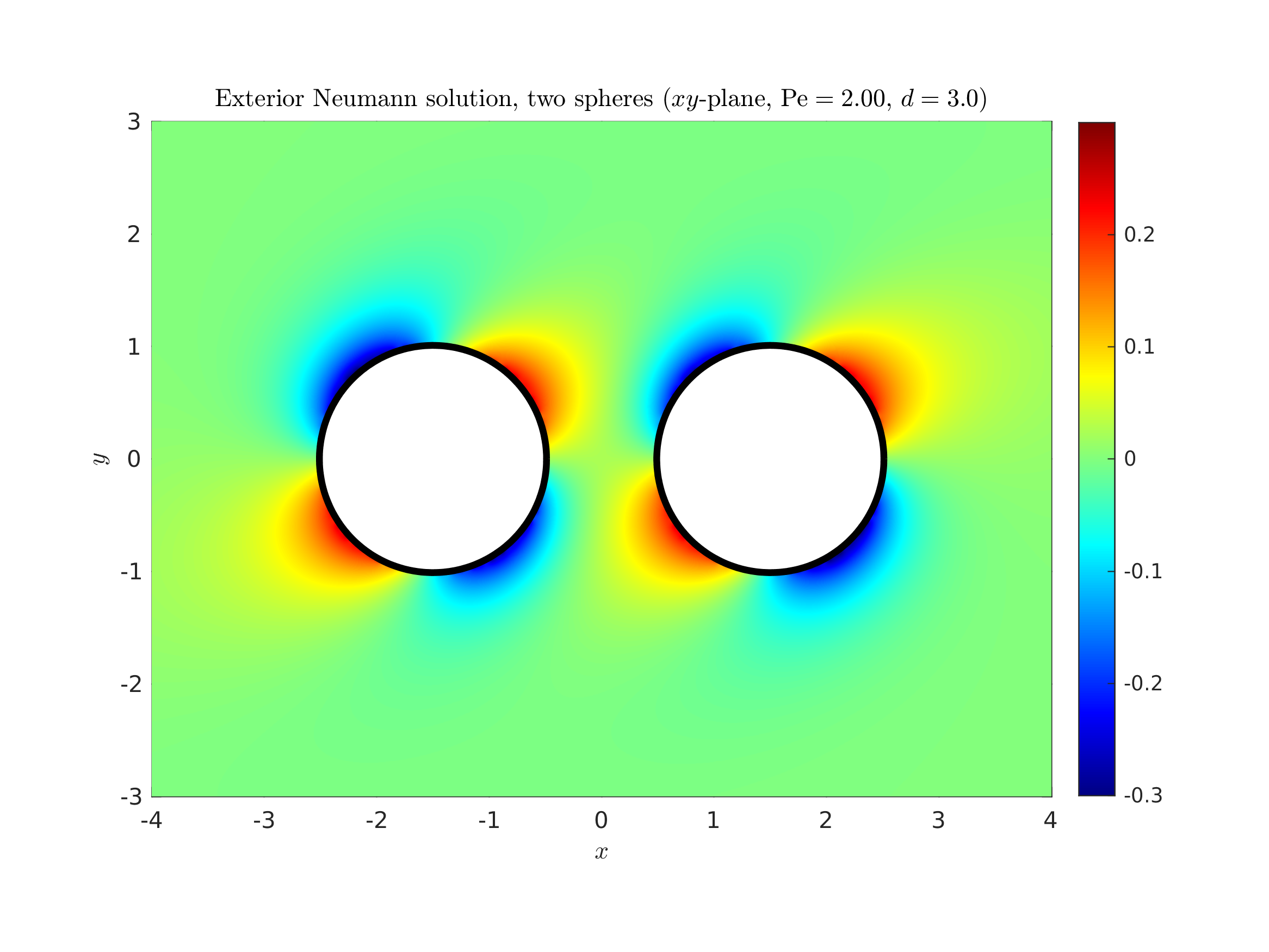}
    \caption{Exterior of two balls.}
    \label{fig:exterior-neumann}
        \end{subfigure}
        \caption{Solution to the exterior Neumann problem from \eqref{eq:smoluchowski_equation_Psi-freq}, $\omega=0$.}
\end{figure}

In this section, we study the exterior Neumann problem \eqref{eq:smoluchowski_equation_Psi-freq} for $\omega = 0$ and different $\Pe$. Two domains $\Omega$ are considered, corresponding to the complement of a ball of radius $1$ centered at $\mathbf{0}$, resp.~two balls of radius $1$ centered in $(\pm \frac{d}{2}, 0, 0)$. 
We solve the exterior Neumann problem \eqref{eq:smoluchowski_equation_Psi-freq} in $\Omega$
using the direct boundary integral formulation \eqref{eq:BIEs_ext} of the second kind: 
$\left( -\tfrac{1}{2}\mathrm{I} + \hat{\mathrm{K}}_0 \right)\varphi = \hat{\mathrm{V}}_0 g$ on $\Gamma$,
with $g(\nex) = -\Pe (\mathbf{n}\cdot \hat{\mathbf{x}}_1)\, x_2$.
Galerkin discretizations in $P_h^1$ are considered,  as in Section \ref{sec:galerkin}. $\Gamma$ is discretized with a sequence of quasi-uniform meshes, as shown in Figure \ref{fig:sphere-mesh} for $3196$ elements. 
Snapshots of the solution in the XY--plane are shown in Figures \ref{fig:exterior-neumann-sphere} (exterior of one ball) and \ref{fig:exterior-neumann} (exterior of two balls) using this mesh. 

This problem arises for colloidal suspensions of hard spherical particles under an externally imposed shear force \cite{blawzdziewicz1993structure} \cite{brady1995normal}. For the study of their rheological behavior  
the components 
\begin{equation}
    Q_{ij}(\Pe) \coloneqq \int\limits_{\Gamma} \hat{\Psi}_0(\nex)x_ix_j\dsx = \int\limits_{\Gamma} \varphi(\nex)x_ix_j\dsx 
\end{equation}
of the stationary shear stress tensor are of key interest.

\begin{figure}[t]
    \centering
    \includegraphics[width=0.4\linewidth]{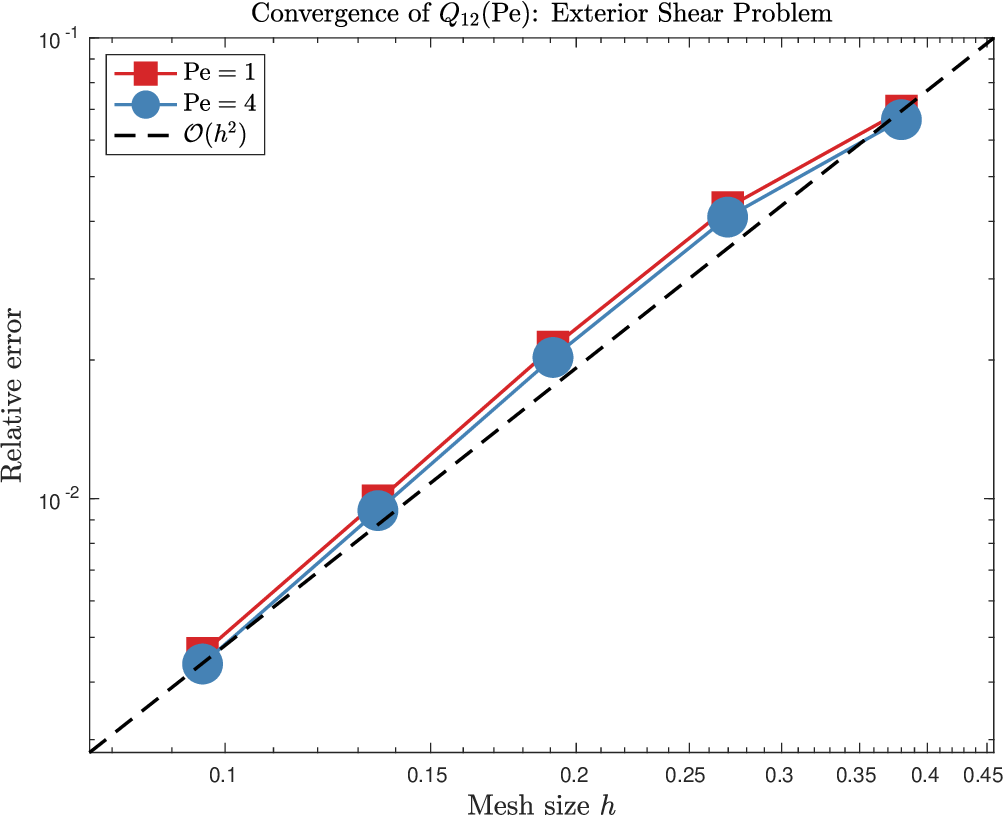}
    \caption{Error in component $Q_{12}$ of shear stress tensor.}
    \label{fig:Q12}
\end{figure}

Figure \ref{fig:Q12} shows the error in the component $Q_{12}$ of the stress tensor as a function of the mesh size $h$. Here, the problem is considered in the geometry with a single ball and $\Pe \in \{1,4\}$, for the stationary case $\omega = 0$. The error in $Q_{12}$ converges quadratically to $0$, as compared to the benchmarks $Q_{12}(\Pe = 1) = 0.26963$, $Q_{12}(\Pe = 4) = 0.93805$ computed for $6396$ elements. 

The methods presented in the article will allow the detailed study of the shear response of colloidal suspensions \cite{future}, beyond the limited accuracy available from finite element approximations \cite{brady1995normal}.  

\section*{Acknowledgments}
We thank Thomas Franosch and Anton L\"{u}ders for discussions about the physical problem and background, as well as Fabian Neu\-hold and Alexander Ostermann for discussions about the numerical analysis.

\bibliographystyle{siamplain}
\bibliography{references}

@article{costabel1988boundary,
  title={Boundary integral operators on {L}ipschitz domains: elementary results},
  author={Costabel, Martin},
  journal={SIAM Journal on Mathematical Analysis},
  volume={19},
  number={3},
  pages={613--626},
  year={1988},
  publisher={SIAM}
}

@article{future,
  title={Exact results for the leading nonanalytic terms in the stress of sheared {B}rownian
hard-sphere suspensions},
  author={Franosch, Thomas and Gimperlein, Heiko and Labarca-Figueroa, Ignacio and L\"{u}ders, Anton and Neuhold, Fabian and Ostermann, Alexander},    journal={in preparation},
  year={2026}
}

@book{bebendorf2008hmatrices,
  author    = {Mario Bebendorf},
  title     = {Hierarchical Matrices: A Means to Efficiently Solve Elliptic Boundary Value Problems},
  series    = {Lecture Notes in Computational Science and Engineering},
  volume    = {63},
  publisher = {Springer},
  address   = {Berlin},
  year      = {2008}
}

@book{hackbusch2015hierarchical,
  author    = {Wolfgang Hackbusch},
  title     = {Hierarchical Matrices: Algorithms and Analysis},
  series    = {Springer Series in Computational Mathematics},
  volume    = {49},
  publisher = {Springer},
  address   = {Heidelberg},
  year      = {2015}
}

@article{gwinnerstephan,
  title={Advanced boundary element methods},
  author={Gwinner, Joachim and Stephan, Ernst Peter},
  journal={Treatment of boundary value, transmission and contact problems. Springer, Cham},
  year={2018},
  publisher={Springer}
}

@incollection{sauter2010boundary,
  title={Boundary element methods},
  author={Sauter, Stefan A and Schwab, Christoph},
  booktitle={Boundary Element Methods},
  pages={183--287},
  year={2010},
  publisher={Springer}
}

@book{gardiner,
  title={Handbook of Stochastic Methods for Physics, Chemistry and the Natural Sciences},
  author={Gardiner, Crispin W.},
  year={1997},
  publisher={Springer}
}

@book{steinbach2008numerical,
  title={Numerical approximation methods for elliptic boundary value problems: finite and boundary elements},
  author={Steinbach, Olaf},
  year={2008},
  publisher={Springer}
}

@article{bergenholtz2002non,
  title={The non-{N}ewtonian rheology of dilute colloidal suspensions},
  author={Bergenholtz, J and Brady, JF and Vicic, M},
  journal={Journal of Fluid Mechanics},
  volume={456},
  pages={239--275},
  year={2002},
  publisher={Cambridge University Press}
}

@article{brady1995normal,
  title={Normal stresses in colloidal dispersions},
  author={Brady, John F and Vicic, Michael},
  journal={Journal of Rheology},
  volume={39},
  number={3},
  pages={545--566},
  year={1995},
  publisher={The Society of Rheology}
}

@article{squires2005simple,
  title={A simple paradigm for active and nonlinear microrheology},
  author={Squires, Todd M and Brady, John F},
  journal={Physics of Fluids},
  volume={17},
  number={7},
  year={2005},
  publisher={AIP Publishing}
}

@book{mclean2000strongly,
  title={Strongly elliptic systems and boundary integral equations},
  author={McLean, William Charles Hector},
  year={2000},
  publisher={Cambridge University Press}
}

@article{franosch2018time,
  title={Time-dependent active microrheology in dilute colloidal suspensions},
  author={Leitmann, Sebastian and Mandal, Suvendu and Fuchs, Matthias and Puertas, Antonio M and Franosch, Thomas},
  journal={Phys. Rev. Fluids},
  volume={3},
  number={10},
  pages={103301},
  year={2018},
  publisher={APS}
}

@article{bruno2017windowed,
  title={Windowed {G}reen function method for the {H}elmholtz equation in the presence of multiply layered media},
  author={Bruno, Oscar P and P{\'e}rez-Arancibia, Carlos},
  journal={Proceedings of the Royal Society A: Mathematical, Physical and Engineering Sciences},
  volume={473},
  number={2202},
  pages={20170161},
  year={2017},
  publisher={The Royal Society Publishing}
}

@article{bruno2016windowed,
  title={Windowed {G}reen function method for layered-media scattering},
  author={Bruno, Oscar P and Lyon, Mark and P{\'e}rez-Arancibia, Carlos and Turc, Catalin},
  journal={SIAM J. Appl. Math.},
  volume={76},
  number={5},
  pages={1871--1898},
  year={2016},
  publisher={SIAM}
}

@article{cheng1998method,
  title={A method of images for the evaluation of electrostatic fields in systems of closely spaced conducting cylinders},
  author={Cheng, Hongwei and Greengard, Leslie},
  journal={SIAM J. Appl. Math.},
  volume={58},
  number={1},
  pages={122--141},
  year={1998},
  publisher={SIAM}
}

@article{gimbutas2015simple,
  title={Simple and efficient representations for the fundamental solutions of {S}tokes flow in a half-space},
  author={Gimbutas, Zydrunas and Greengard, Leslie and Veerapaneni, Shravan},
  journal={J. Fluid Mech.},
  volume={776},
  pages={R1},
  year={2015},
  publisher={Cambridge University Press}
}

@article{bruno2014rapidly,
  title={Rapidly convergent two-dimensional quasi-periodic Green function throughout the spectrum—including Wood anomalies},
  author={Bruno, Oscar P and Delourme, B{\'e}rang{\`e}re},
  journal={Journal of Computational Physics},
  volume={262},
  pages={262--290},
  year={2014},
  publisher={Elsevier}
}

@article{lin2024greensII,
  title={The {G}reen's function for an acoustic, half-space impedance problem, {P}art {II}: Analysis of the slowly varying and the plane wave component},
  author={Lin, C and Melenk, Jens Markus and Sauter, S},
  journal={arXiv:2408.03587},
  year={2024}
}

@article{lin2024greensI,
  title={An Explicit Factorization of the {G}reen's Function for an Acoustic Half-Space Problem with Impedance Boundary Conditions into an Oscillatory Exponential and a Slowly Varying Function},
  author={Lin, C and Melenk, Jens Markus and Sauter, S},
  journal={Math. Methods Appl. Sci.},
  volume={48},
  number={13},
  pages={12807--12812},
  year={2025},
  publisher={Wiley Online Library}
}

@article{lindbo2011spectral,
  title={Spectral accuracy in fast {E}wald-based methods for particle simulations},
  author={Lindbo, Dag and Tornberg, Anna-Karin},
  journal={J. Comp. Phys.},
  volume={230},
  number={24},
  pages={8744--8761},
  year={2011},
  publisher={Elsevier}
}

@article{arens2013analysing,
  title={Analysing {E}wald's method for the evaluation of {G}reen's functions for periodic media},
  author={Arens, Tilo and Sandfort, Kai and Schmitt, Susanne and Lechleiter, Armin},
  journal={IMA J. Appl. Math.},
  volume={78},
  number={3},
  pages={405--431},
  year={2013},
  publisher={Oxford University Press}
}

@article{liang2025accelerating,
  title={Accelerating Fast {E}wald Summation with Prolates for Molecular Dynamics Simulations},
  author={Liang, Jiuyang and Lu, Libin and Barnett, Alex and Greengard, Leslie and Jiang, Shidong},
  journal={arXiv:2505.09727},
  year={2025}
}

@article{broms2025method,
  title={A method of fundamental solutions for large-scale 3D elastance and mobility problems},
  author={Broms, Anna and Barnett, Alex H and Tornberg, Anna-Karin},
  journal={Adv. Comp. Math.},
  volume={51},
  number={5},
  pages={45},
  year={2025},
  publisher={Springer}
}

@article{babuvska1994p,
  title={The p and h-p versions of the finite element method, basic principles and properties},
  author={Babu{\v{s}}ka, Ivo and Suri, Manil},
  journal={SIAM Review},
  volume={36},
  number={4},
  pages={578--632},
  year={1994},
  publisher={SIAM}
}

@book {schwab1998hpfem,
    AUTHOR = {Schwab, Ch.},
     TITLE = {{$p$}- and {$hp$}-finite element methods},
    SERIES = {Numerical Mathematics and Scientific Computation},
 PUBLISHER = {The Clarendon Press, Oxford University Press, New York},
      YEAR = {1998},
     PAGES = {xii+374},
      ISBN = {0-19-850390-3}
}

@article {zhang2021hpquad,
    AUTHOR = {Zhang, Mingzhu and Mao, Xinyu and Yi, Lijun},
     TITLE = {Exponential convergence of the {$hp$}-version of the composite
              {G}auss-{L}egendre quadrature for integrals with endpoint
              singularities},
   JOURNAL = {Appl. Numer. Math.},
  FJOURNAL = {Applied Numerical Mathematics. An IMACS Journal},
    VOLUME = {170},
      YEAR = {2021},
     PAGES = {340--352}
     }

@article{elrick1962source,
  title={Source functions for diffusion in uniform shear flow},
  author={Elrick, DE},
  journal={Australian Journal of Physics},
  volume={15},
  pages={283},
  year={1962}
}

@book{dhont1996introduction,
  title={An introduction to dynamics of colloids},
  author={Dhont, Jan KG},
  year={1996},
  publisher={Elsevier}
}

@article{brader2010nonlinear,
  title={Nonlinear rheology of colloidal dispersions},
  author={Brader, Joseph M},
  journal={Journal of Physics: Condensed Matter},
  volume={22},
  number={36},
  pages={363101},
  year={2010}
}

@article{blawzdziewicz1993structure,
  title={Structure and rheology of semidilute suspension under shear},
  author={B{\l}awzdziewicz, J and Szamel, Grzegorz},
  journal={Physical Review E},
  volume={48},
  number={6},
  pages={4632},
  year={1993},
  publisher={APS}
}
\end{document}